\documentclass[article,10pt]{elsarticle}
\usepackage{amssymb}
\usepackage{amsmath}
\usepackage{amsfonts}
\usepackage{amsthm}
\usepackage{graphics,graphicx}
\usepackage{float}
\usepackage{mathtools}
\usepackage[dvipsnames]{xcolor}
\usepackage{url}

\usepackage[english]{babel}

\journal{ADAM}
\newcommand{\qbinom}[2]{\genfrac{[}{]}{0pt}{}{#1}{#2}}
\newtheorem{theorem}{Theorem}
\newtheorem{lemma}{Lemma}

\newtheorem{corollary}{Corollary}
\newtheorem{example}{Example}
\newtheorem{conjecture}{Conjecture}
\newtheorem{remark}{Remark}
\newtheorem{problem}{Problem}
\newtheorem{o-problem}{Open problem}
\begin{document}

\begin{frontmatter}

\title{Minimum supports of eigenfunctions of graphs: a survey \tnoteref{label0}}
\tnotetext[label0]{The study was carried out within the framework of the state contract
of the Sobolev Institute of Mathematics (project no. 0314-2019-0016).}

\author[01]{Ev Sotnikova}
\ead{ev.v.sotnikova@gmail.com}

\author[01]{Alexandr Valyuzhenich\corref{cor1}}
\cortext[cor1]{Corresponding author}
\ead{graphkiper@mail.ru}

\address[01]{Sobolev Institute of Mathematics, Ak. Koptyug av. 4, Novosibirsk 630090, Russia}

\begin{abstract}
In this work we present a survey of results on the problem of finding the minimum cardinality of the support of eigenfunctions of graphs.
\end{abstract}

\begin{keyword}
eigenfunction\sep eigenfunctions of graphs\sep eigenspace\sep minimum support\sep trade\sep bitrade\sep $1$-perfect bitrade\sep  weight distribution bound
\vspace{\baselineskip}
\MSC[2010] 05C50\sep 05E30\sep 05B30\sep 15A18
\end{keyword}

\end{frontmatter}
\section{Introduction}\label{SectionIntro}
The eigenvalues of a graph are closely related to its structural properties and invariants (see the monographs \cite{BrouwerHaemers,CvetkovicDS,CvetkovicRS10}).
Eigenfunctions (equivalently, eigenvectors) of graphs, in contrast to their eigenvalues, have received only sporadic attention of researchers.
In particular, basic properties of eigenfunctions of graphs can be found in the work of Merris \cite{Merris}. Among the most famous results we can recall the theory around Perron-Frobenius vector \cite{Per07, Fro12, BRN97} with its applications to a variety of problems including ranking, population growth models, Markov chains behavior and many other \cite{PageBrin99, Keener93, Mac00, Seneta06}; and the results about Fiedler vector \cite{F75, Chung97, Das04} and its connection to the problems of spectral graph partitioning and clustering \cite{PSL90, NJW01}, graph coloring \cite{AG84}, graph drawing \cite{KCH03} and other (for example, \cite{SM00, SBJ07}).
In addition, it is worth noting a series of works \cite{Biyikoglu03,Biyikoglu04,Biyikoglu05,Colin,Davies,Duval,Friedman,Powers,vanHolst}
devoted to various discrete versions of Courant's nodal domain theorem. We refer the reader to \cite{BiyikogluBook}, \cite[Chapter 9]{CvetkovicRS97} and \cite{Saito} for more details about eigenfunctions of graphs.

In this work we consider undirected graphs without loops and multiple edges.
The eigenvalues of a graph are the eigenvalues of its adjacency matrix.
Let $G=(V,E)$ be a graph with vertex set $V=\{v_1,\ldots,v_n\}$ and let $\lambda$ be an eigenvalue of $G$. The set of neighbors of a vertex $x$ is denoted by $N(x)$.
A function $f:V\longrightarrow{\mathbb{R}}$ is called a {\em $\lambda$-eigenfunction} of $G$ if $f\not\equiv 0$ and the equality
\begin{equation}\label{EqEigenfunction}
\lambda\cdot f(x)=\sum_{y\in{N(x)}}f(y)
\end{equation}
holds for any vertex $x\in V$.
Note that if $f$ is a $\lambda$-eigenfunction of $G$, then $A\overrightarrow{f}=\lambda \overrightarrow{f}$, where
$A$ is the adjacency matrix of $G$ and 
$\overrightarrow{f}=(f(v_1),\ldots,f(v_n))^T$, i.e. $\overrightarrow{f}$ is an eigenvector of the matrix $A$
with eigenvalue $\lambda$.
The set of functions $f:V\longrightarrow{\mathbb{R}}$ satisfying (\ref{EqEigenfunction}) for any vertex $x\in V$ is called a {\em $\lambda$-eigenspace} of $G$. Denote by $U_{\lambda}(G)$ the $\lambda$-eigenspace of $G$.
The {\em support} of a function $f:V\longrightarrow{\mathbb{R}}$ is the set $S(f)=\{x\in V~|~f(x)\neq 0\}$.
A $\lambda$-eigenfunction of $G$ is called {\em optimal} if it has the minimum cardinality of the support among all $\lambda$-eigenfunctions of $G$.
In this work we focus on the following extremal problem for eigenfunctions of graphs.
\begin{problem}[MS-problem]\label{ProblemMS}
Let $G$ be a graph and let $\lambda$ be an eigenvalue of $G$. Find the minimum cardinality of the support of a $\lambda$-eigenfunction of $G$.
\end{problem}
In what follows, in this work we will use the abbreviation MS-problem instead of Problem \ref{ProblemMS}.
Now we discuss the deep connection between MS-problem and the intersection problem of
two combinatorial objects and the problem of finding the minimum size of trades.

Many combinatorial objects (equitable partitions, completely regular codes, Steiner systems $S(k-1,k,n)$, $1$-perfect codes, etc.) can be defined as eigenfunctions of graphs with some discrete restrictions.
The study of such objects often leads to the problem of finding the minimum possible difference between two objects from the same class
(for example, see \cite{EtzionVardy,FranklPach,Hwang,Potapov10,PotapovPPI}).
Since the symmetric difference of such two objects is also an eigenfunction of the corresponding graph, this problem is directly related to MS-problem.

Trades of different types are used for constructing and studying the structure of different combinatorial objects (combinatorial $t$-designs, codes, Latin squares, etc.).
Trades are also studied independently as some natural generalization of objects of the corresponding type
(trades can exist even if the corresponding complete objects do not exist).
Roughly speaking, trades reflect possible differences between two combinatorial objects from the same class:
if $C'$ and $C''$ are two combinatorial objects with the same parameters, then the pair $(C'\setminus C'',C''\setminus C')$ is a trade
(for more information on trades see \cite{Billington,Cavenagh,HedayatKhosrovshahi,KrotovMogilnykhPotapov}).
Many types of trades ($T(k-1,k,v)$ Steiner trades, $q$-ary $T_q(k-1,k,v)$ Steiner trades, $1$-perfect trades, extended $1$-perfect trades, latin trades, etc.) can be represented as eigenfunctions of
the corresponding graphs with some additional discrete restrictions (for example, see \cite[Section 2.4]{KrotovPerfectBitrades}).
So, for such trades the problem of finding the minimum size can be reduced to MS-problem for the corresponding graphs (see, for example, \cite{KrotovMogilnykhPotapov,Valyuzhenich20,VorobevKrotov}).

In particular, MS-problem has appeared as a natural generalization of the following results.
\begin{itemize}
  \item Let $C_1$ and $C_2$ be two distinct binary perfect codes of length $n=2^m-1$.
  In \cite{EtzionVardy} Etzion and Vardy proved that the maximum possible cardinality of their intersection $C_1\cap C_2$ is $2^{n-m}-2^{\frac{n-1}{2}}$.
  Equivalently, they found the minimum possible cardinality of their symmetric difference $C_1\triangle C_2$.
This result can be proved by applying the so-called weight distribution bound for the Hamming graph $H(n,2)$ and its eigenvalue $-1$ (see Subsection \ref{SubSectionEquitable} and Section~\ref{sec:WDB}).
  \item In \cite{Hwang} Hwang proved that the minimum size of a $T(t,k,v)$ trade is $2^{t+1}$ and obtained a characterization
  of $T(t,k,v)$ trades of size $2^{t+1}$.
  In particular, the minimum size of a $T(t,k,v)$ Steiner trade was found in \cite{Hwang}. For $t=k-1$ this result can be proved by applying the weight distribution bound for the Johnson graph $J(v,k)$ and its eigenvalue $-k$ (see Subsection \ref{SubSectionSteinerTrades} and Section \ref{sec:WDB}). It is interesting that Frankl and Pach \cite{FranklPach} also found  the minimum size of a $T(t,k,v)$ trade. They formulated their results in terms of null $t$-designs. In Section~\ref{sec:grassmann} we will meet null designs again during our discussion about optimal eigenfunctions of the Grassmann graph.

\end{itemize}

MS-problem was first formulated by Krotov and Vorob'ev \cite{VorobevKrotov} in 2014 (they considered MS-problem for the Hamming graph).
During the last six years, MS-problem has been actively studied for various families of distance-regular graphs
\cite{Bespalov,GoryainovKabanovShalaginovValyuzhenich,Krotovtezic,KrotovMogilnykhPotapov,SotnikovaCubical,SotnikovaBilinear,Valyuzhenich,Valyuzhenich20,ValyuzhenichVorobev,VorobevKrotov,VMV} and Cayley graphs on the symmetric group \cite{KabanovKonstantinovaShalaginovValyuzhenich}.
In particular, MS-problem is completely solved for all eigenvalues of the Hamming graph \cite{Valyuzhenich20,ValyuzhenichVorobev} and asymptotically solved for all eigenvalues of the Johnson graph \cite{VMV}. Note that for eigenfunctions of distance-regular graphs a lower bound for its support cardinality is known. This bound is called the weight distribution bound and we will discuss it in details in Section~\ref{sec:WDB}.
In this work we give a survey of results on MS-problem. We also discuss constructions of optimal eigenfunctions and 
the main ideas of the proofs of the results.

Now we would like to consider the following problem.
\begin{problem}\label{ProblemMMS}
Let $G=(V,E)$ be a graph and let $\lambda$ be an eigenvalue of $G$. Find
$$\min_{f\in U_{\lambda}(G),f\not\equiv 0}|\{x\in V~|~f(x)\geq 0\}|.$$
\end{problem}

Note, that the statements of MS-problem and Problem~\ref{ProblemMMS} are similar. An analogue of Problem \ref{ProblemMMS} for association schemes was first formulated in 1984 by Bier \cite{Bier84}.
Later, Bier and Delsarte \cite{BierDelsarte1,BierDelsarte2} and Bier \cite{Bier88} studied the same problem for eigenvectors belonging to the
direct sum of several eigenspaces of an association scheme.
Bier and Manickam \cite{BierManickam}, Manickam and Mikl\'{o}s \cite{ManickamMiklos} and Manickam and Singhi \cite{ManickamSinghi} initiated the study of Problem \ref{ProblemMMS} for the second largest eigenvalue of Johnson and Grassmann graphs.
In particular, the following two conjectures were formulated in 1988.

\begin{conjecture}[Manickam, Mikl\'{o}s and Singhi \cite{ManickamMiklos,ManickamSinghi}]\label{ConjectureMMS}
Let $x_1,\ldots,x_n$ be real numbers such that $x_1+\ldots+x_n=0$.
If $n\geq 4k$, then there are at least $\binom{n-1}{k-1}$
$k$-element subsets of the set $\{x_1,\ldots,x_n\}$ with nonnegative sum.
\end{conjecture}

The second conjecture is an analogue of Conjecture \ref{ConjectureMMS} for vector spaces.
Let $V$ be an $n$-dimensional vector space over a finite field $\mathbb{F}_q$.
Let $\qbinom{V}{k}_q$ denote the family of all $k$-dimensional subspaces of $V$ and let
$\qbinom{n}{k}_q$ denote the $q$-Gaussian binomial coefficient.
For each $1$-dimensional subspace $v\in \qbinom{V}{1}_q$, assign a real-valued weight $f(v)\in \mathbb{R}$ so that
the sum of all weights is zero. For a general subspace $S\subset V$, define its weight $f(S)$ to be
the sum of the weights of all the $1$-dimensional subspaces it contains.

\begin{conjecture}[Manickam and Singhi \cite{ManickamSinghi}]\label{ConjectureVectorMMS}
Let $V$ be an $n$-dimensional vector space over $\mathbb{F}_q$ and let $f:\qbinom{V}{1}_q\rightarrow \mathbb{R}$  be a weighting of the $1$-dimensional subspaces such that $\sum_{v\in \qbinom{V}{1}_q}f(v)=0$.
If $n\geq 4k$, then there are at least $\qbinom{n-1}{k-1}_q$ $k$-dimensional subspaces with nonnegative weight.
\end{conjecture}

Conjecture \ref{ConjectureMMS} is still open. However, there are several relatively recent works \cite{Alon,Chowdhury,Frankl,Pokrovskiy} with polynomial bounds. In particular, Alon, Huang and Sudakov \cite{Alon} verified Conjecture \ref{ConjectureMMS} for $n\geq 33k^2$.
A linear bound $n\geq 10^{46}k$ was obtained by Pokrovskiy \cite{Pokrovskiy}.
In 2014 Chowdhury, Sarkis and Shahriari \cite{Chowdhury} and Huang and Sudakov \cite{Huang} independently showed that Conjecture \ref{ConjectureVectorMMS} holds for $n\geq 3k$. Using the technique of the work \cite{Chowdhury}, Ihringer \cite{Ihringer} proved that Conjecture \ref{ConjectureVectorMMS} is true for $n\geq 2k$ and large $q$.
Some new results on Problem \ref{ProblemMMS} for the third largest eigenvalue of the Johnson graph can be found in \cite{MVV}.
It seems very intriguing to establish the interconnection between Problem \ref{ProblemMMS} and MS-problem.

The paper is organized as follows. In Section \ref{SectionEigenfunctionsInConfig}, we give two examples of combinatorial problems that are closely related to MS-problem. In Section \ref{SectionDefinitions}, we introduce basic definitions and notations.
In Section \ref{sec:WDB}, we discuss what the weight distribution bound is and how it can be calculated from the intersection arrays of the distance-regular graphs. We complete this section with several intuitive examples. In Sections \ref{SectionHamming}-\ref{SectionStar}, we give a survey of results on MS-problem for the Hamming graph, the Doob graph, the Johnson graph, the Grassmann graph, the bilinear forms graph, the Paley graph and the Star graph respectively. In Section \ref{SectionRemarks}, we present some observations on optimal eigenfunctions of graphs. 
In Section \ref{SectionOpenProblems}, we formulate several open problems.
\section{Eigenfunctions in combinatorial configurations and MS-problem}\label{SectionEigenfunctionsInConfig}
In this section, we recall that equitable $2$-partitions, $1$-perfect codes and $T(k-1,k,v)$ Steiner trades can be defined as eigenfunctions of graphs with some discrete restrictions. We also discuss the connections of MS-problem with the intersection problem of
two $1$-perfect codes of a given graph and the problem of finding the minimum size of Steiner trades.
\subsection{Equitable partitions and $1$-perfect codes}\label{SubSectionEquitable}
Let $G=(V,E)$ be a graph.
An ordered $r$-partition $(C_1,\ldots, C_{r})$ of $V$ is called {\em equitable} if for any
$i,j\in \{1,\ldots,r\}$ there is $S_{i,j}$ such that any vertex of $C_i$
has exactly $S_{i,j}$ neighbors in $C_j$.
The matrix $S=(S_{i,j})_{i,j\in
\{1,\ldots,r\}}$ is called the {\em quotient matrix} of the
equitable partition. A set $C\subseteq V$ is called a {\em $1$-perfect code} in $G$ if every ball of radius $1$ contains one vertex from $C$.
For more information on equitable partitions and perfect codes we refer the reader to \cite{BespalovKMTVPerfectColor}, \cite[Chapter 5]{Godsil} and
\cite{Ahlswede,Heden,Solov'eva08,Solov'eva13}.

Let $G$ be a $k$-regular graph and let $(C_1,C_2)$ be an equitable $2$-partition of $G$ with the quotient matrix
$$S=\begin{pmatrix}
a & b\\
c & d\\
\end{pmatrix}.$$
The eigenvalues of $S$ are $k$ and $a-c$.
We define the function $f_{(C_1,C_2)}$ on the vertices of $G$ by the following rule:
$$
f_{(C_1,C_2)}(x)=\begin{cases}
b,&\text{if $x\in C_1$;}\\
-c,&\text{if $x\in C_2$.}
\end{cases}
$$
One can verify that $f_{(C_1,C_2)}$ is an $(a-c)$-eigenfunction of $G$. So, any equitable $2$-partition can be represented as an eigenfunction of the corresponding graph.
Suppose that $C$ is a $1$-perfect code in $G$. Then the partition $(C,\overline{C})$ is equitable with the quotient matrix
$$\begin{pmatrix}
0 & k\\
1 & k-1\\
\end{pmatrix}.$$
Therefore, the function $f_{(C,\overline{C})}$ is a $(-1)$-eigenfunction of $G$. So, if $C_1$ and $C_2$ are $1$-perfect codes in $G$,
then the function $f=f_{(C_1,\overline{C_1})}-f_{(C_2,\overline{C_2})}$ is also a $(-1)$-eigenfunction of $G$. Moreover, we have the equality
$$|S(f)|=|C_1\triangle C_2|.$$
Thus, the problem of finding the minimum cardinality of the symmetric difference of two distinct $1$-perfect codes of a regular graph can be reduced to MS-problem for this graph and eigenvalue $-1$.
\subsection{$T(k-1,k,v)$ Steiner trades}\label{SubSectionSteinerTrades}
Let $v$, $k$, $t$ be positive integers such that $v>k>t$ and let $X$ be a set of size $v$.
A pair $(T_0,T_1)$ of disjoint collections of $k$-subsets (blocks) of $X$ is called a {\em $T(t,k,v)$ trade} if every $t$-subset of $X$ is included in
the same number of blocks of $T_0$ and $T_1$. The {\em size} of a $T(t,k,v)$ trade $(T_0,T_1)$ is $|T_0|+|T_1|$. A $T(t,k,v)$ trade is called {\em Steiner} if every $t$-subset of $X$ is included in at most one block of
$T_0$ ($T_1$). For further details on $T(t,k,v)$ trades we refer the reader to \cite{Billington,HedayatKhosrovshahi,KhosrovshahiMT}.

Suppose that $(T_0,T_1)$ is a $T(k-1,k,v)$ Steiner trade. The {\em Johnson graph} $J(v,k)$ can be defined as follows.
The vertices of $J(v,k)$ are $k$-subsets of $X$, and two vertices are adjacent if they have exactly $k-1$ common elements.
We define the function $f_{(T_0,T_1)}$ on the vertices of $J(v,k)$ by the following rule:
$$
f_{(T_0,T_1)}(x)=\begin{cases}
1,&\text{if $x\in T_0$;}\\
-1,&\text{if $x\in T_1$;}\\
0,&\text{otherwise.}
\end{cases}
$$
For a $(k-1)$-subset $A$ of $X$ denote by $C(A)$ the set of vertices of $J(v,k)$ containing the set $A$ (these vertices form a clique of size $v-k+1$ in $J(v,k)$).
We note that $C(A)$ either contains one element from $T_0$ and one element from $T_1$ or does not contain elements from $T_0\cup T_1$. Using this fact, one can easily check that $f_{(T_0,T_1)}$ is a $(-k)$-eigenfunction of $J(v,k)$.
Moreover, we have the equality $$|S(f_{(T_0,T_1)})|=|T_0|+|T_1|.$$
Thus, the problem of finding the minimum size of $T(k-1,k,v)$ Steiner trades can be reduced to MS-problem
for the Johnson graph $J(v,k)$ and its eigenvalue $-k$.

\section{Basic definitions}\label{SectionDefinitions}

Recall that a distance $d_G(v,\,u)=d(u,\,v)$ between two vertices $v$ and $u$ in a graph $G=(V,\,E)$ is the length of the shortest path that connects them. The largest distance between any pairs of vertices is called the diameter $D$. A connected graph $G=(V,\,E)$ is called distance-regular if it is regular of degree $k$ and for any two vertices $v,\,u\in V$ at distance $i=d(v,\,u)$ there are precisely $c_i$ neighbors of $u$ which are at distance $i-1$ from $v$ and precisely $b_i$ neighbors of $u$ which are at distance $i+1$ from $v$; where $c_i$ and $b_i$ do not depend on the choice of vertices $u$ and $v$ but depend only on $d(u,\,v)$. Numbers $b_i$, $c_i$, $a_i=k-b_i-c_i$ are called the intersection numbers and a set $\{b_0,b_1,\ldots,b_{D-1};\,c_1,\ldots,c_D\}$ is called an intersection array of a distance-regular graph $G$. For more details about distance-regular graphs, the reader is referred to a classical monograph \cite{BCN89} and a recent survey \cite{DKT16}.

Let $G_1=(V_1,E_1)$ and $G_2=(V_2,E_2)$ be simple graphs.
The {\em Cartesian product} $G_1\square G_2$ of graphs $G_1$ and $G_2$ is defined as follows.
The vertex set of $G_1\square G_2$ is $V_1\times V_2$;
and any two vertices $(x_1,y_1)$ and $(x_2,y_2)$ are adjacent if and only if either
$x_1=x_2$ and $y_1$ is adjacent to $y_2$ in $G_2$, or
$y_1=y_2$ and $x_1$ is adjacent to $x_2$ in $G_1$.

Suppose $G_1=(V_1,E_1)$ and $G_2=(V_2,E_2)$ are two graphs. Let $f_1:V_1\longrightarrow{\mathbb{R}}$ and $f_2:V_2\longrightarrow{\mathbb{R}}$.
Denote $G=G_1\square G_2$.
We define the {\em tensor product} $f_1\cdot f_2$  on the vertices of $G$ by the following rule:
$$(f_1\cdot f_2)(x,y)=f_1(x)f_2(y)$$ for $(x,y)\in V(G)=V_1\times V_2$.
We will use the tensor product of functions for constructing optimal eigenfunctions of the Hamming and Doob graphs in Subsection \ref{SubSectionConsrHamming} and Section \ref{SectionDoob}.

Let $\rm{Sym}(X)$ denote the symmetric group on a finite set $X$ and let $\rm{Sym}_{n}$ denote the symmetric group on the set $\{1,\ldots,n\}$.

Let $\Sigma_q=\{0,1,\ldots,q-1\}$. Let $f(x_1,\ldots,x_n)$ be a function defined on the set $\Sigma_{q}^{n}$, let $\pi\in \rm{Sym}_n$ and let $\sigma_1,\ldots,\sigma_{n}\in \rm{Sym}(\Sigma_q)$.  We define the functions $f_{\pi}$ and $f_{\pi,\sigma_{1},\ldots,\sigma_{n}}$ as follows:
$$f_{\pi}(x_1,\ldots,x_n)=f(x_{\pi(1)},\ldots,x_{\pi(n)})$$
and
$$f_{\pi,\sigma_{1},\ldots,\sigma_{n}}(x_1,\ldots,x_n)=f(\sigma_{1}(x_{\pi(1)}),\ldots,\sigma_{n}(x_{\pi(n)})).$$
We will use the functions $f_{\pi}$ and $f_{\pi,\sigma_{1},\ldots,\sigma_{n}}$ in Subsections \ref{SubSectionConsrHamming} and \ref{SubSectionResults}.

Let $G=(V,E)$ be a graph.
A set $C\subseteq V$ is called a {\em completely regular code} in $G$ if the partition
$(C^{(0)},\ldots,C^{(\rho)})$ is equitable, where $C^{(d)}$ is the set of vertices at distance $d$ from $C$ and $\rho$ ({\em the covering radius} of $C$) is the maximum $d$ for which $C^{(d)}$ is nonempty.
In other words, a subset of $V$ is a completely regular code in $G$
if the distance partition with respect to the subset is equitable.
For more information on completely regular codes see \cite{BorgesRifaZinoviev}, \cite[Chapter 11.7]{Godsil} and \cite{KoolenKM,KoolenLM}.
We will use completely regular codes in Section \ref{SectionRemarks}.
\section{Weight Distribution Bound}\label{sec:WDB}
In this section we recall what a weight distribution bound is and how it can be used as a lower bound for MS-problem in case of distance-regular graphs. 

Weight distribution bound is well known and has appeared in several papers under different disguise (for more details see \cite{K11}, \cite{KrotovMogilnykhPotapov}). In order for this survey to be self-contained we would like to provide the full proof and equip the reader with several intuitive examples.


Let $A$ be the adjacency matrix of some distance-regular graph $G=(V,\,E)$.  Now consider the distance-$t$ graph $G_t=(V,\,E_t)$ defined as follows: two vertices $v$ and $u$ are adjacent in $G_t$ if and only if they are at distance $t$ in $G$. In other words, $\{v,\, u\}\in E_t \iff d_G(v,\, u) = t$. By $A_t$ we denote the adjacency matrix of $G_t$.

Considering the combinatorial definition of distance regularity from the matrix point of view, we obtain the following recurrence (see, for example, equation (1) in \cite{DKT16}): 
\begin{equation}\label{eq:DistRegRec}
    A_tA=a_tA_t+b_{t-1}A_{t-1}+c_{t+1}A_{t+1},
\end{equation}

for $t=0,1,\ldots, D$ where $b_{-1}A_{-1}=c_{D+1}A_{D+1}=0$. 

From the above we can show that there exist polynomials $P_t$ of degree $t$ such that:
\[A_t=P_t(A),\quad t=0,1,\ldots,D.\]
	
It is well known that in case of Hamming graphs these polynomials are actually Kravchuk polynomials (up to some linear change of variables) and those are Eberlein polynomials in case of Johnson graphs.

But how can we make use of it in finding the lower bound for our MS-problem? Suppose $f$ is a $\lambda$-eigenfunction of our graph $G$. Since $A^i \overrightarrow{f}=\lambda^i \overrightarrow{f}$, we get the following equations:
\[A_t\overrightarrow{f}=P_t(A)\overrightarrow{f}=P_t(\lambda) \overrightarrow{f}.\]

In other words, $f$ is a $P_t(\lambda)$-eigenfunction of graph $G_t$. As an immediate consequence we obtain:

\[P_t(\lambda)f(v)=\sum\limits_{\substack{u\in V,\\d(u,\,v)=t}}f(u).\]
	
In other words, in distance-regular graphs the sum of the eigenfunction values on the vertices at distance $t$ from a fixed vertex $v$ depends only on $f(v)$ and the corresponding eigenvalue. Without lost of generality we can consider $f(v)=1$. The array $[1,P_1(\lambda),\ldots,P_D(\lambda)]$ is called the weight distribution of a $\lambda$-eigenfunction. 
	
Thus from (\ref{eq:DistRegRec}) we can write the following recurrence:
\[
\begin{gathered}
P_0(\lambda)=1,\\
P_1(\lambda)=\lambda,\\
P_t(\lambda)=\frac{\lambda P_{t-1}(\lambda)-b_{t-2}P_{t-2}(\lambda)-a_{t-1}P_{t-1}(\lambda)}{c_t},\mbox{ where }t=2,\ldots,D.
\end{gathered}
\]

Now we are just one step away from obtaining the lower bound we are looking for. Let $w$ be such a vertex that $|f(w)|=\underset{u\in V}{\max} |f(u)|$. The trick is the following: instead of $\lambda$-eigenfunction $f$ we consider a function

\[
g = \frac{1}{f(w)}f.
\]

Thus $g$ is also a $\lambda$-eigenfunction and $S(f)=S(g)$. Moreover, $g(w)=1$ and $|g(u)|\le g(w)=1$ for all $u\in V$. From the weight distribution we obtain \[P_t(\lambda)=\sum\limits_{d(u,\,w)=t}g(u).\]

Therefore, $g$ has at least $|P_t(\lambda)|$ non-zero values at distance $t$ from a vertex $w$. This proves the next lemma.

\begin{lemma}[\cite{KrotovMogilnykhPotapov}, Corollary~1]
	Let $f$ be a $\lambda$-eigenfunction for a distance-regular graph $G$ of diameter $D$, then the following bound takes place:
	\[|S(f)|\ge \sum_{i=0}^D |P_i(\lambda)|.\]
\end{lemma}

In case of irrational eigenvalues this bound can be refined:

\begin{lemma}
	Let $f$ be a $\lambda$-eigenfunction for a distance-regular graph $G$ of diameter $D$, then the following bound takes place:
	\[|S(f)|\ge \sum_{i=0}^D\lceil |P_i(\lambda)|\rceil.\]
\end{lemma}

Let us illustrate this technique on some well-known graphs (see \cite{SotnikovaCubical} for details). We start with the Petersen graph. The Petersen graph is a cubical distance-regular graph on $10$ vertices. Its intersection array is $\{3,2; 1,1\}$ and its eigenvalues are $\{-2^{(4)}, 1^{(5)}, 3^{(1)}\}$. Calculating the weight distribution we obtain $[1, \lambda, \lambda^2-3]$.
\begin{itemize}
    \item For $\lambda=1$ it gives us the lower bound $4$. An optimal $1$-eigenfunction achieves this bound. A subgraph induced on non-zero vertices can be described as two non-incident edges. An example is presented below (Figure~\ref{fig:Petersen1}).
    
    \begin{figure}[H]
    \begin{center}
    \includegraphics[scale=0.4]{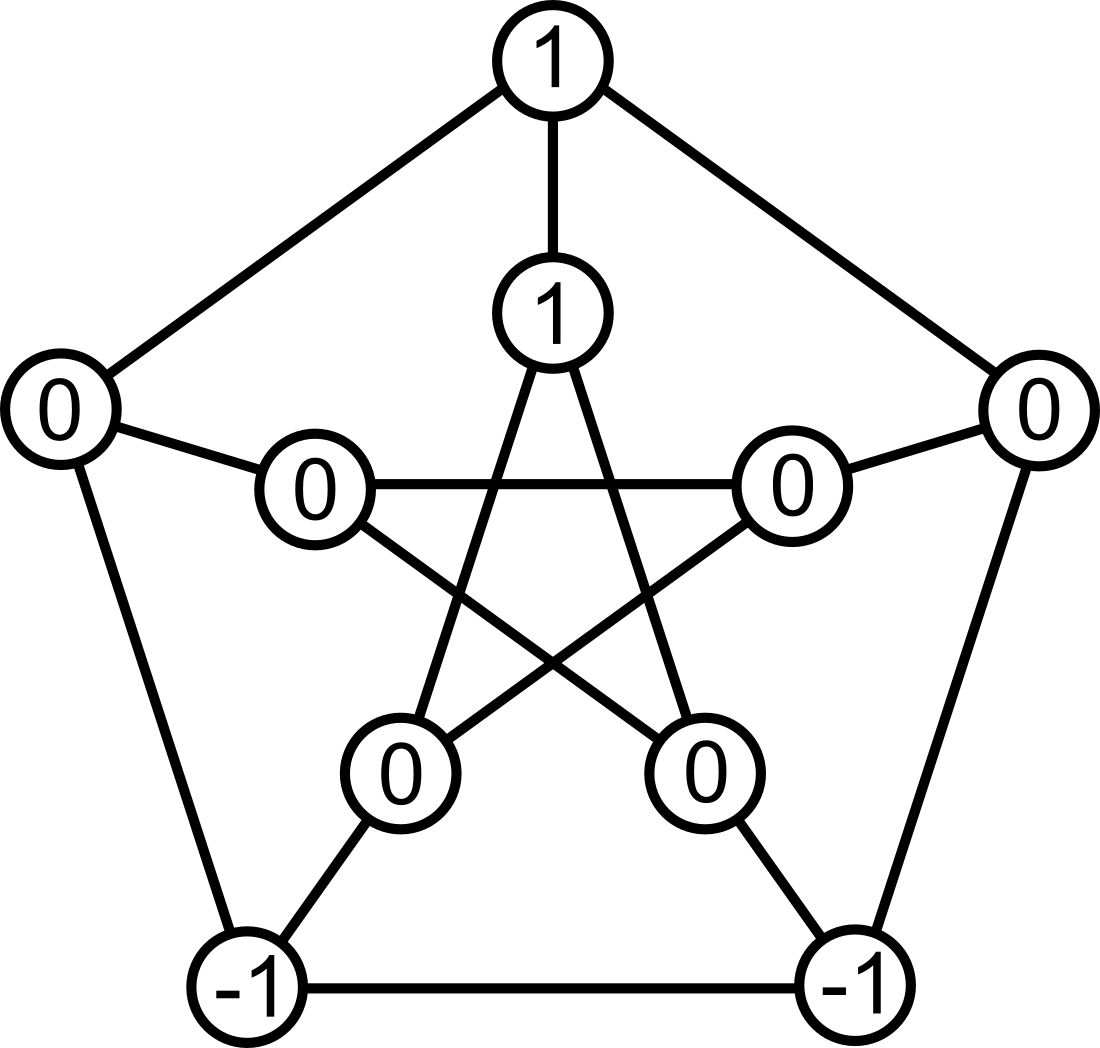}
    \end{center}
    \caption{Optimal $1$-eigenfunction of the Petersen graph.}\label{fig:Petersen1}
    \end{figure}
    
    \item For $\lambda=-2$ the lower bound is the same. But this case is different because this bound cannot be achieved. Optimal $(-2)$-eigenfunction has a support of cardinality $6$ and the corresponding induced subgraph is either a cycle on six vertices, or $H$-graph. See Figure~\ref{fig:petersen-cycle} and Figure~\ref{fig:petersen-hgraph}.
    
    \begin{figure}[H]
    \begin{minipage}[h]{0.45\linewidth}
    \center{\includegraphics[width=0.75\linewidth]{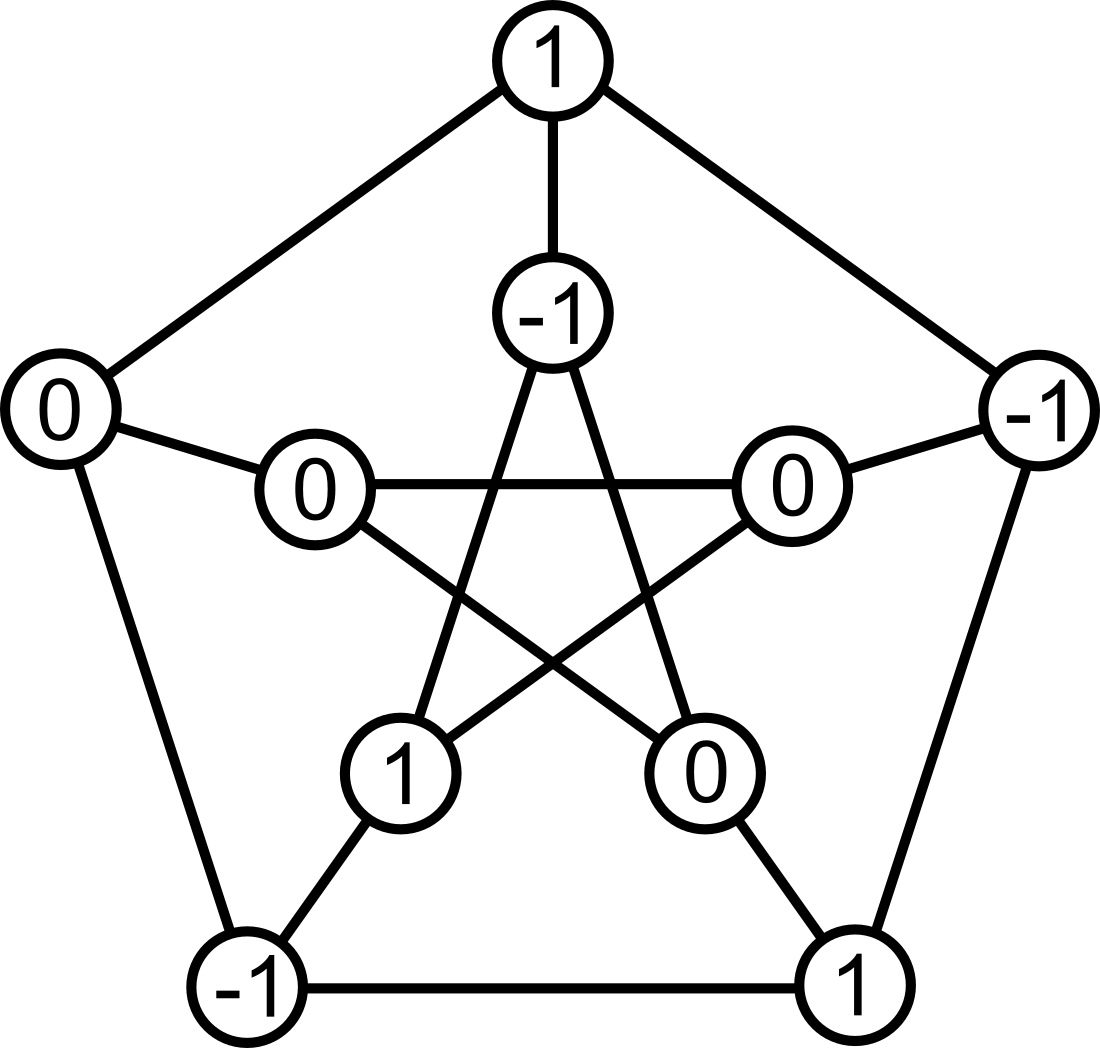}}
    \caption{Optimal $(-2)$-eigenfunction of the Petersen graph --- cycle}
    \label{fig:petersen-cycle}
    \end{minipage}
    \hfill
    \begin{minipage}[h]{0.45\linewidth}
    \center{\includegraphics[width=0.75\linewidth]{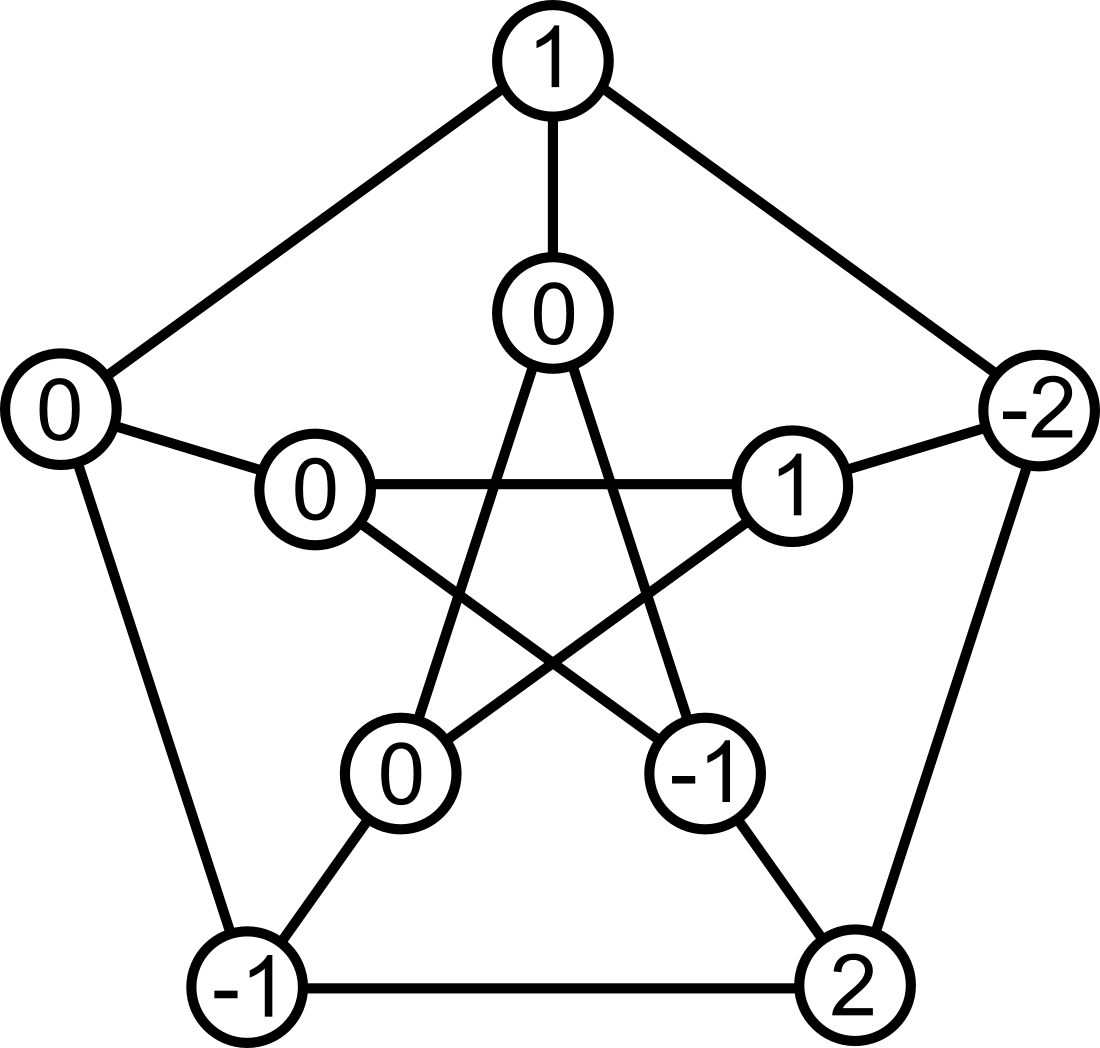}}
    \caption{Optimal $(-2)$-eigenfunction of the Petersen graph --- H-graph}
    \label{fig:petersen-hgraph}
    \end{minipage}
    \end{figure}
    \end{itemize}

\begin{figure}[H]
    \begin{center}
    \includegraphics[scale=0.3]{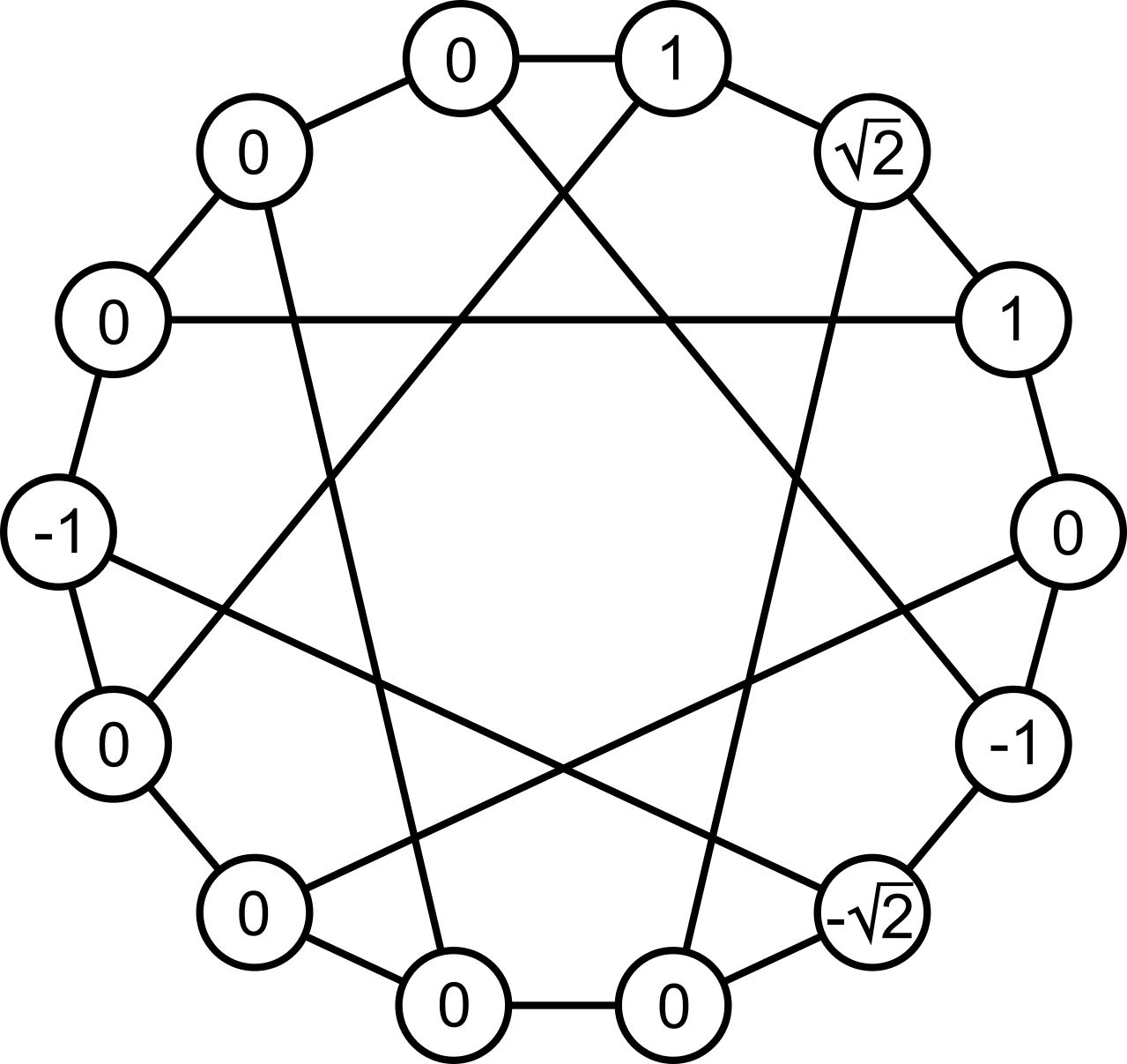}
    \end{center}
    \caption{Optimal $\sqrt2$-eigenfunction of the Heawood graph.}\label{fig:heawood}
    \end{figure}

As a quick illustration of bound refinement, let us consider the Heawood graph, a distance-regular graph on $14$ vertices. Its intersection array is $\{3,2,2;1,1,3\}$ and its spectrum is $\{\pm3^{(1)}, \pm\sqrt{2}^{(6)}\}$. The weight distribution is $[1,\lambda, \lambda^2-3, \frac{1}{3}(\lambda^3-5\lambda)]$. Thus for $\lambda=\pm\sqrt 2$ the exact weight distribution bound is $2+2\sqrt2$, while the refined bound is $6$. Fig.~\ref{fig:heawood} presents an example of an optimal $\sqrt{2}$-eigenfunction.

More examples can be found in \cite{SotnikovaCubical}, where MS-problem is solved together with a characterisation of such functions for $10$ out of $13$ cubical distance-regular graphs for all their eigenvalues.

Thus for any distance-regular graph a lower bound on a cardinality of a $\lambda$-eigenfunction support can be calculated directly from the intersection array of a graph with respect to the corresponding eigenvalue $\lambda$. However this bound is not necessary feasible. We will see in the next sections that a weight distribution bound is achieved for:
\begin{itemize}
    \item an eigenvalue $-1$ of the boolean Hamming graph of an odd order and a minimum eigenvalue of arbitrary Hamming graphs;
    \item a minimum eigenvalue of the Paley graphs of square order (since this graph is self-complementary of diameter 2, this property also holds for a second non-principal eigenvalue);
    \item a minimum eigenvalue of the Johnson graphs;
    \item a minimum eigenvalue of the Grassmann graphs;
    \item a minimum eigenvalue of strongly regular bilinear forms graphs over a prime field.
\end{itemize}
However this is not the case for a minimum eigenvalue of bilinear forms graphs of larger diameter, for example.

We can see from the list above that the minimum eigenvalue $\lambda_D$ attracted some special attention. This is because we can say something about non-zeros of optimal $\lambda_D$-eigenfunctions. Indeed for every distance-regular graph admitting a so-called Delsarte pair the existence of a $\lambda_D$-eigenfunction $f$ achieving the weight distribution bound is equivalent to the existence of an isometric distance-regular subgraph induced on the non-zero values of $f$ (see Corollary~2 from \cite{KrotovMogilnykhPotapov}). We will explore this property more in the Section~\ref{sec:bfg}, but actually this result is a partial case of a more general theory which also demonstrates the strong connection between eigenfunctions and bitrades. For a deeper dive into the topic the interested reader is referred to Sections~2 and 3 of \cite{KrotovMogilnykhPotapov}.

\section{Hamming graph}\label{SectionHamming}
In this section, we give a survey of results on MS-problem and its generalizations for the Hamming graph.
The {\em Hamming graph} $H(n,q)$ is defined as follows. Let $\Sigma_q=\{0,1,\ldots,q-1\}$.
The vertex set of $H(n,q)$ is $\Sigma_{q}^n$, and two vertices are adjacent if they differ in exactly one position.
This graph is a distance-regular graph.
The Hamming graph $H(n,q)$ has $n+1$ distinct eigenvalues $\lambda_i(n,q)=n(q-1)-q\cdot i$, where $0\leq i\leq n$.
Denote by $U_{i}(n,q)$ the $\lambda_{i}(n,q)$-eigenspace of $H(n,q)$. The direct sum of subspaces
$$U_i(n,q)\oplus U_{i+1}(n,q)\oplus\ldots\oplus U_j(n,q)$$ for $0\leq i\leq j\leq n$ is denoted by $U_{[i,j]}(n,q)$.
We say that a function $f\in U_{[i,j]}(n,q)$, where $f\not\equiv 0$, is {\em optimal} in the space $U_{[i,j]}(n,q)$ if
$|S(f)|\leq |S(g)|$ for any function $g\in U_{[i,j]}(n,q)$, $g\not\equiv 0$.

Firstly, we briefly discuss all results on MS-problem for the Hamming graph. After that, we will consider the more general Problem \ref{GMSHamming} for functions from the space $U_{[i,j]}(n,q)$.
In \cite{Krotovtezic} Krotov based on the approach of work \cite{Potapov} proved that the minimum cardinality of the support of a $\lambda_{i}(n,2)$-eigenfunction of $H(n,2)$ is $\max (2^i,2^{n-i})$.
In \cite{VorobevKrotov} Krotov and Vorob'ev showed that the cardinality of the support of a $\lambda_{i}(n,q)$-eigenfunction of
$H(n,q)$ is at least $$2^i\cdot (q-2)^{n-i}$$ for $\frac{iq^2}{2n(q-1)}>2$ and
$$q^n\cdot (\frac{1}{q-1})^{i/2}\cdot (\frac{i}{n-i})^{i/2}\cdot (1-\frac{i}{n})^{n/2}$$ for $\frac{iq^2}{2n(q-1)}\leq2$.
In \cite{Valyuzhenich} Valyuzhenich for $q\geq 3$ proved that the minimum cardinality of the support of a $\lambda_{1}(n,q)$-eigenfunction of $H(n,q)$ is
$2\cdot (q-1)\cdot q^{n-2}$ and obtained a characterization of optimal $\lambda_{1}(n,q)$-eigenfunctions.
Later in \cite{Valyuzhenich20,ValyuzhenichVorobev} the following generalization of MS-problem for the Hamming graph was considered.
\begin{problem}\label{GMSHamming}
Let $n\geq 1$, $q\ge 2$ and $0\leq i\leq j\leq n$. Find the minimum cardinality of the support of functions from the space $U_{[i,j]}(n,q)$.
\end{problem}
In \cite{ValyuzhenichVorobev} Valyuzhenich and Vorob'ev found the minimum cardinality of the support of a function from the space $U_{[i,j]}(n,q)$
for arbitrary $q\geq 3$ except the case when $q=3$ and $i+j>n$.
Moreover, in \cite{ValyuzhenichVorobev} a characterization of functions that are optimal in the space $U_{[i,j]}(n,q)$
was obtained for $q\ge 3$, $i+j\le n$ and $q\ge 5$, $i=j$, $i>\frac{n}{2}$.
In \cite{Valyuzhenich20} Valyuzhenich found the minimum cardinality of the support of a function from the space $U_{[i,j]}(n,q)$ for $q=2$ and $q=3$, $i+j>n$. Thus, Problem \ref{GMSHamming} is completely solved for all $n\geq 1$ and $q\ge 2$. As a consequence, MS-problem for the Hamming graph is also solved for all eigenvalues.

In what follows, in this section we will consider in detail Problem \ref{GMSHamming}.
In Subsection \ref{SubSectionConsrHamming}, we present constructions of functions that are optimal in the space $U_{[i,j]}(n,q)$.
In Subsection \ref{SubSectionResults}, we give a survey of results on Problem \ref{GMSHamming} and discuss the main ideas of the proof of these results. In particular, we carefully explore Lemma \ref{ReductionHamming} which is a key tool for solving Problem \ref{GMSHamming}.
In Subsection \ref{SubSectionTrade}, we focus on a connection between Problem \ref{GMSHamming} and the problem of finding the minimum size of $1$-perfect bitrades in the Hamming graph.
\subsection{Constructions of functions with the minimum cardinality of the support}\label{SubSectionConsrHamming}
In this subsection, we discuss constructions of functions that are optimal in the space $U_{[i,j]}(n,q)$.
It is interesting that in all cases such functions are constructed as a tensor product of several elementary optimal functions defined on the vertices of the Hamming graph of diameter not greater than three.

Firstly, we define five sets of elementary optimal functions.

\begin{figure}[H]
\begin{center}
\includegraphics[scale=0.37]{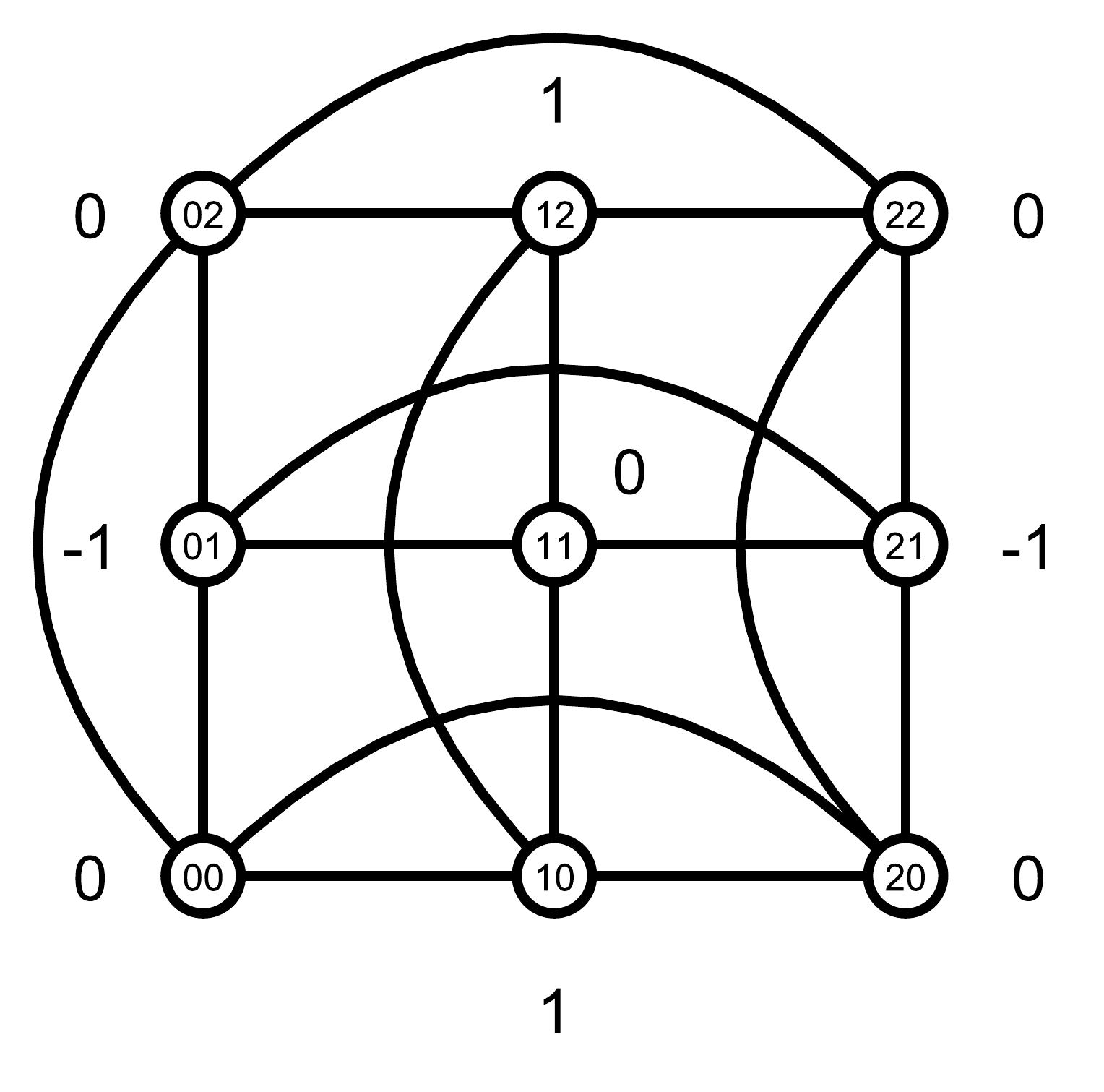}
\end{center}
\caption{Function $a_{3,1,1}$ in $H(2,3)$.}\label{Figure1}
\end{figure}

For $k,m\in{\Sigma_q}$ we define the function $a_{q,k,m}$ on the vertices of the Hamming graph $H(2,q)$ by the following rule:
$$
a_{q,k,m}(x,y)=\begin{cases}
1,&\text{if $x=k$ and $y\neq m$;}\\
-1,&\text{if $y=m$ and $x\neq k$;}\\
0,&\text{otherwise.}
\end{cases}
$$
The function $a_{3,1,1}$ is shown in Figure \ref{Figure1}. We note that $a_{q,k,m}$ is optimal in the space $U_{1}(2,q)$ for any $k,m\in{\Sigma_q}$. Denote $A_q=\{a_{q,k,m}~|~k,m\in{\Sigma_q}\}$.

We define the function $\varphi_1$ on the vertices of the Hamming graph $H(2,3)$ by the following rule:
$$
\varphi_1(x,y)=\begin{cases}
1,&\text{if $x=y=0$;}\\
-1,&\text{if $x=1$ and $y=2$;}\\
0,&\text{otherwise.}
\end{cases}
$$

\begin{figure}[H]
\begin{center}
\includegraphics[scale=0.27]{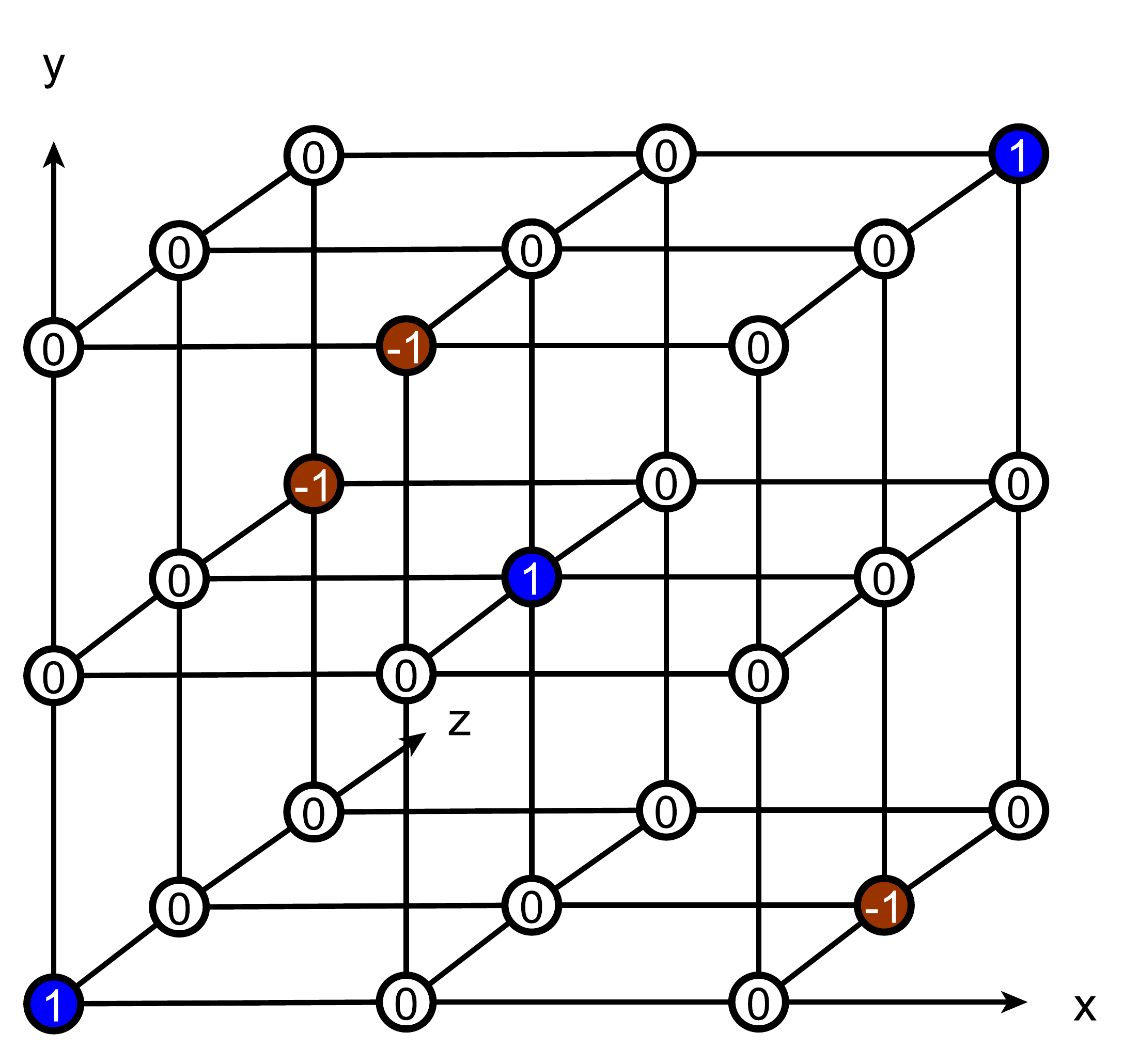}
\end{center}
\caption{Function $\varphi(x,y,z)$ in $H(3,3)$.}\label{Figure2}
\end{figure}

For $a,b\in{\Sigma_{3}}$ denote by $a\oplus b$ the sum of $a$ and $b$ modulo $3$.
We define the function $\varphi$ on the vertices of the Hamming graph $H(3,3)$ by the following rule:
$$
\varphi(x,y,z)=\begin{cases}
\varphi_1(x,y),&\text{if $z=0$;}\\
\varphi_1(x\oplus1,y\oplus1),&\text{if $z=1$;}\\
\varphi_1(x\oplus2,y\oplus2),&\text{if $z=2$.}
\end{cases}
$$
The function $\varphi$ is shown in Figure \ref{Figure2}. We note that  $\varphi$ is optimal in the space $U_2(3,3)$.
Denote $$B=\{\varphi_{\pi,\sigma_{1},\sigma_{2},\sigma_{3}}~|~\pi\in \rm{Sym}_3,\sigma_{1},\sigma_{2},\sigma_{3}\in \rm{Sym}(\Sigma_3)\}.$$

\begin{figure}[H]
\begin{center}
\includegraphics[scale=0.33]{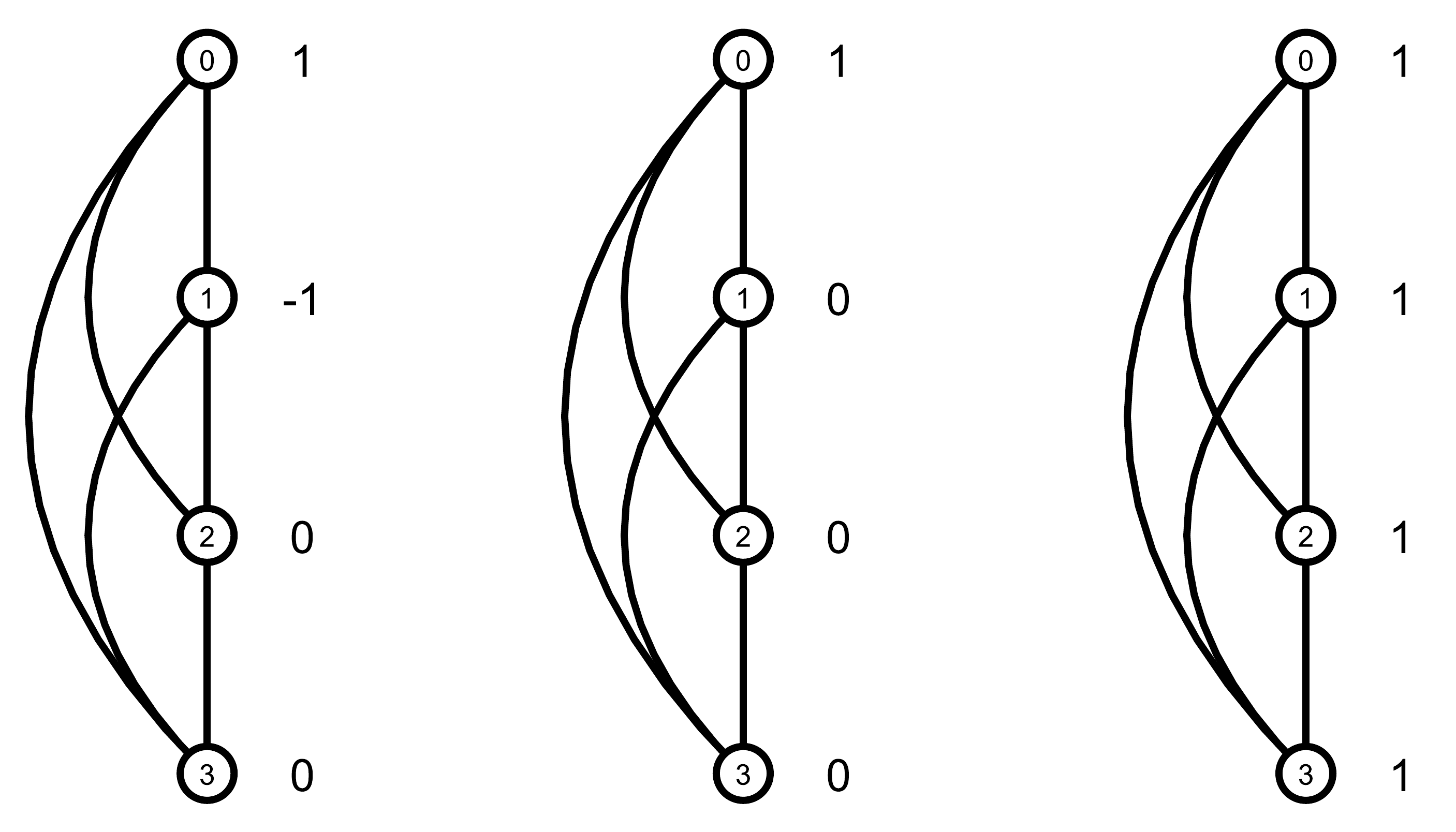}
\end{center}
\caption{Functions $c_{4,0,1}$, $d_{4,0}$ and $e_4$ in $H(1,4)$.}\label{Figure3}
\end{figure}

For $k,m\in{\Sigma_q}$ and $k\neq m$ we define the function $c_{q,k,m}$ on the vertices of the Hamming graph $H(1,q)$ by the following rule:
$$
c_{q,k,m}(x)=\begin{cases}
1,&\text{if $x=k$;}\\
-1,&\text{if $x=m$;}\\
0,&\text{otherwise.}
\end{cases}
$$
The function $c_{4,0,1}$ is shown in Figure \ref{Figure3}. We note that $c_{q,k,m}$ is optimal in the space $U_{1}(1,q)$ for any $k,m\in{\Sigma_q}$ and $k\neq m$.
Denote $C_q=\{c_{q,k,m}~|~k,m\in{\Sigma_q},k\neq m\}$.

For $k\in{\Sigma_q}$ we define the function $d_{q,k}$ on the vertices of the Hamming graph $H(1,q)$ by the following rule:
$$
d_{q,k}(x)=\begin{cases}
1,&\text{if $x=k$;}\\
0,&\text{otherwise.}
\end{cases}
$$
The function $d_{4,0}$ is shown in Figure \ref{Figure3}. We note that $d_{q,k}$ is optimal in the space $U_{[0,1]}(1,q)$ for any $k\in{\Sigma_q}$.
Denote $D_q=\{d_{q,k}~|~k\in{\Sigma_q}\}$.

Let $e_q:\Sigma_{q}\longrightarrow{\mathbb{R}}$ and $e_q\equiv 1$. The function $e_{4}$ is shown in Figure \ref{Figure3}.
We note that $e_q$ is optimal in the space $U_{0}(1,q)$.
Denote $E_q=\{e_q\}$.

Now, we define four classes of functions that are optimal in the space $U_{[i,j]}(n,q)$ for the corresponding cases.

Let $i+j\le n$. We say that a function $f$ defined on the vertices of $H(n,q)$ belongs to the class $\mathcal{F}_{1}(n,q,i,j)$ if $$f=c\cdot \prod_{k=1}^{i}g_{k}\cdot \prod_{k=1}^{n-i-j}h_{k}\cdot \prod_{k=1}^{j-i}v_{k},$$ where $c$ is a real non-zero constant, $g_k\in{A_q}$ for $k\in{[1,i]}$, $h_k\in{E_q}$ for $k\in{[1,n-i-j]}$ and $v_k\in{D_q}$ for $k\in{[1,j-i]}$.

Let $i+j>n$. We say that a function $f$ defined on the vertices of $H(n,q)$ belongs to the class $\mathcal{F}_{2}(n,q,i,j)$ if $$f=c\cdot \prod_{k=1}^{n-j}g_{k}\cdot \prod_{k=1}^{i+j-n}h_{k}\cdot \prod_{k=1}^{j-i}v_{k},$$ where $c$ is a real non-zero constant, $g_k\in{A_q}$ for $k\in{[1,n-j]}$, $h_k\in{C_q}$ for $k\in{[1,i+j-n]}$ and $v_k\in{D_q}$ for $k\in{[1,j-i]}$.

Let $\frac{i}{2}+j\leq n$ and $i+j>n$. We say that a function $f$ defined on the vertices of $H(n,3)$ belongs to the class $\mathcal{F}_{3}(n,i,j)$ if $$f=c\cdot \prod_{k=1}^{2n-i-2j}g_{k}\cdot \prod_{k=1}^{i+j-n}h_{k}\cdot \prod_{k=1}^{j-i}v_{k},$$ where $c$ is a real non-zero constant, $g_k\in{A_3}$ for $k\in{[1,2n-i-2j]}$, $h_k\in{B}$ for $k\in{[1,i+j-n]}$ and $v_k\in{D_3}$ for $k\in{[1,j-i]}$.

Let $\frac{i}{2}+j>n$. We say that a function $f$ defined on the vertices of $H(n,3)$ belongs to the class $\mathcal{F}_{4}(n,i,j)$ if $$f=c\cdot \prod_{k=1}^{n-j}g_{k}\cdot \prod_{k=1}^{i+2j-2n}h_{k}\cdot \prod_{k=1}^{j-i}v_{k},$$ where $c$ is a real non-zero constant, $g_k\in{B}$ for $k\in{[1,n-j]}$, $h_k\in{C_3}$ for $k\in{[1,i+2j-2n]}$ and $v_k\in{D_3}$ for $k\in{[1,j-i]}$.

We note that functions from $\mathcal{F}_{1}(n,q,i,j)$ and $\mathcal{F}_{2}(n,q,i,j)$ are optimal in the space $U_{[i,j]}(n,q)$ for $q\geq 2$, $i+j\leq n$ and $q\geq 2$ ($q\neq 3$), $i+j>n$ respectively.
We also note that functions from $\mathcal{F}_{3}(n,i,j)$ and $\mathcal{F}_{4}(n,i,j)$ are optimal in the space $U_{[i,j]}(n,3)$ for
$\frac{i}{2}+j\leq n$, $i+j>n$ and $\frac{i}{2}+j>n$ respectively.

\subsection{Problem \ref{GMSHamming}}\label{SubSectionResults}
In this subsection, we discuss Problem \ref{GMSHamming}.
The following theorem is a combination of the results proved in \cite{Valyuzhenich20,ValyuzhenichVorobev} (see \cite[Theorems 1 and 3]{ValyuzhenichVorobev} and \cite[Theorems 3-6]{Valyuzhenich20}).

\begin{theorem}\label{TheoremHamming}

\begin{enumerate}
  \item Let $f\in{U_{[i,j]}(n,q)}$, where $q\geq 2$, $i+j\le n$ and $f\not\equiv 0$. Then $$|S(f)|\geq 2^{i}\cdot(q-1)^{i}\cdot q^{n-i-j}$$ and this bound is sharp. Moreover, for $q\geq 3$ the equality

      $|S(f)|=2^{i}\cdot (q-1)^{i}\cdot q^{n-i-j}$ holds if and only if $f_{\pi}\in \mathcal{F}_{1}(n,q,i,j)$ for some permutation $\pi\in \rm{Sym}_n$.

  \item Let $f\in{U_{[i,j]}(n,q)}$, where $q\geq 2$, $q\neq 3$, $i+j>n$ and $f\not\equiv 0$. Then $$|S(f)|\geq 2^{i}\cdot(q-1)^{n-j}$$ and this bound is sharp. Moreover, for $i=j$ and $q\geq 5$ the equality $|S(f)|=2^{i}\cdot(q-1)^{n-i}$ holds if and only if $f_{\pi}\in \mathcal{F}_{2}(n,q,i,i)$ for some permutation $\pi\in \rm{Sym}_n$.

  \item Let $f\in{U_{[i,j]}(n,3)}$, where $\frac{i}{2}+j\leq n$, $i+j>n$ and $f\not\equiv 0$. Then $$|S(f)|\geq 2^{3(n-j)-i}\cdot3^{i+j-n}$$ and this bound is sharp.

  \item Let $f\in{U_{[i,j]}(n,3)}$, where $\frac{i}{2}+j> n$ and $f\not\equiv 0$. Then $$|S(f)|\geq 2^{i+j-n}\cdot3^{n-j}$$ and this bound is sharp.
\end{enumerate}
\end{theorem}
Now, we discuss the main ideas of the proof of Theorem \ref{TheoremHamming}.

Let $f$ be a real-valued function defined on the vertices of the Hamming graph $H(n,q)$ and let $k\in{\Sigma_q}$, $r\in\{1,\ldots,n\}$.
We define a function $f_{k}^{r}$ on the vertices of $H(n-1,q)$ as follows:
for any vertex $y=(y_1,\ldots,y_{r-1},y_{r+1},\ldots,y_n)$ of $H(n-1,q)$
$$f_{k}^{r}(y)=f(y_1,\ldots,y_{r-1},k,y_{r+1},\ldots,y_n).$$
One of the important points in the proof of Theorem \ref{TheoremHamming} is the following.
\begin{lemma}[\cite{ValyuzhenichVorobev}, Lemma 4]\label{ReductionHamming}
Let $f\in{U_{[i,j]}(n,q)}$ and $r\in\{1,2,\ldots,n\}$. Then the following statements are true:
\begin{enumerate}
\item $f_{k}^{r}-f_{m}^{r}\in{U_{[i-1,j-1]}(n-1,q)}$ for $k,m\in\Sigma_q$.
\item $\sum_{k=0}^{q-1}f_{k}^{r}\in{U_{[i,j]}(n-1,q)}$.
\item $f_{k}^{r}\in{U_{[i-1,j]}(n-1,q)}$ for $k\in\Sigma_q$.
\end{enumerate}
\end{lemma}
Lemma \ref{ReductionHamming} is a very useful tool for studying of eigenfunctions of the Hamming graph.
It shows the connection between eigenspaces of the Hamming graphs $H(n,q)$ and $H(n-1,q)$.
In particular, this lemma allows to apply induction on $n$, $i$ and $j$ (we can use the induction assumption for the functions $f_{k}^{r}-f_{m}^{r}$, $\sum_{k=0}^{q-1}f_{k}^{r}$ and $f_{k}^{r}$).
Moreover, we suppose that Lemma \ref{ReductionHamming} can be useful not only for the MS-problem but also for other problems.
For example, recently in \cite{MV} Mogilnykh and Valyuzhenich used Lemma \ref{ReductionHamming} for investigation of equitable $2$-partitions of the Hamming graph with the eigenvalue $\lambda_2(n,q)$.
One interesting generalization of Lemma \ref{ReductionHamming} for the products of graphs can be found in \cite[Theorem 3.11]{Taranenko}.
\subsection{Minimum 1-perfect bitrades in the Hamming graph}\label{SubSectionTrade}
In this subsection, we discuss one interesting application of Theorem \ref{TheoremHamming} for the problem of finding the minimum size of $1$-perfect bitrades in the Hamming graph.

Let us recall some definitions. Let $G=(V,E)$ be a graph. For a vertex $x\in V$ denote $B(x)=N(x)\cup \{x\}$.
Let $T_0$ and $T_1$ be two disjoint nonempty subsets of $V$.
The ordered pair $(T_0,T_1)$ is called a {\em $1$-perfect bitrade} in $G$ if for any vertex $x\in V$ the set $B(x)$ either contains one element from $T_0$ and one element from $T_1$ or does not contain elements from $T_0\cup T_1$. The {\em size} of a $1$-perfect bitrade $(T_0,T_1)$ is $|T_0|+|T_1|$.

\begin{figure}[H]
\begin{center}
\includegraphics[scale=0.4]{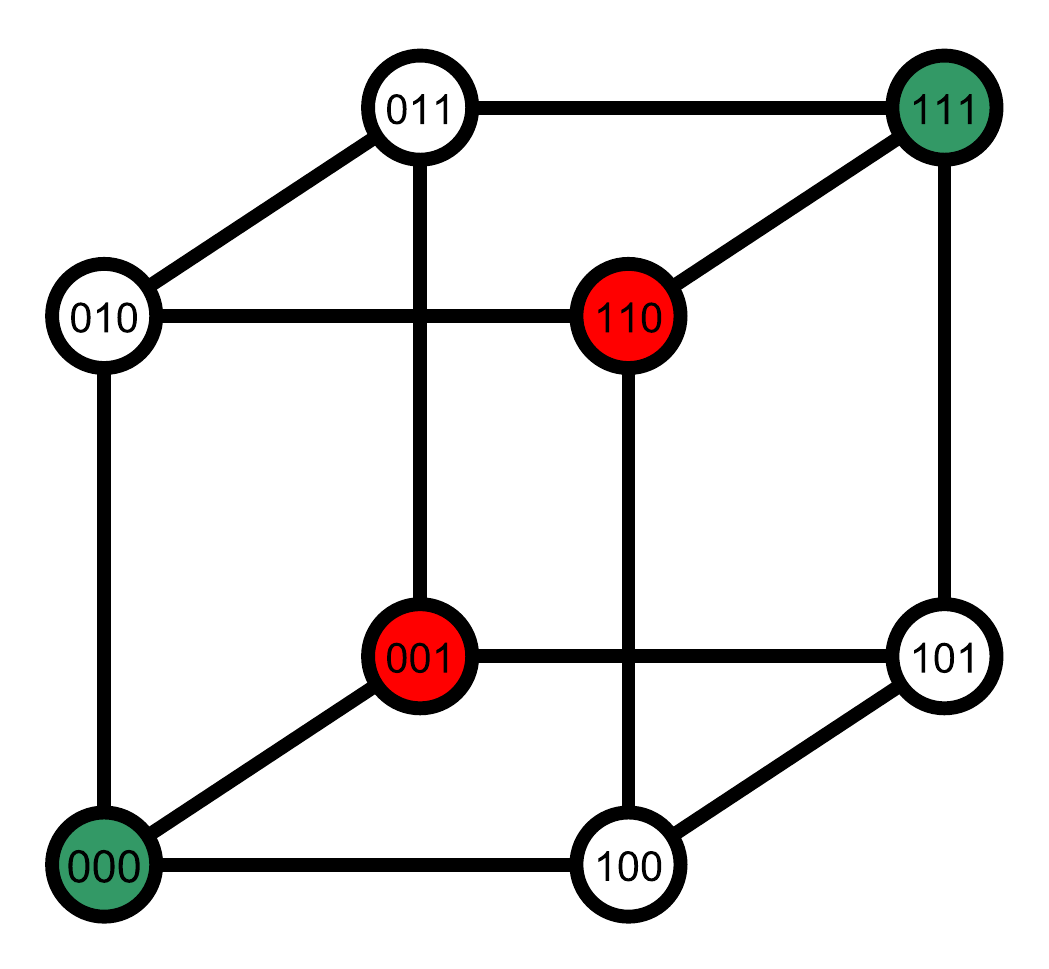}
\end{center}
\caption{$1$-perfect bitrade in $H(3,2)$.}\label{Figure4}
\end{figure}

\begin{example}\label{ExamplePerfectBitrade}
Let $T_0=\{000,111\}$ and $T_1=\{001,110\}$. Then $(T_0,T_1)$ is a $1$-perfect bitrade of size $4$ in $H(3,2)$ (see Figure \ref{Figure4}).
\end{example}

\begin{example}\label{ExamplePerfectBitradeCodes}
Let $G=(V,E)$ be a graph. Suppose $C_1$ and $C_2$ be two distinct $1$-perfect codes in $G$. Then $(C_1\setminus C_2,C_2\setminus C_1)$ is a $1$-perfect bitrade in $G$.
\end{example}

Let $(T_0,T_1)$ be a $1$-perfect bitrade in a graph $G=(V,E)$. We define the function $f_{(T_0,T_1)}:V\longrightarrow{\{-1,0,1\}}$ by the following rule:
$$
f_{(T_0,T_1)}(x)=\begin{cases}
1,&\text{if $x\in T_0$;}\\
-1,&\text{if $x\in T_1$;}\\
0,&\text{otherwise.}
\end{cases}
$$
In what follows, in this subsection we will consider the following problem.
\begin{problem}\label{ProblemPerfectBitrade}
Let $n\ge 3$ and $q\ge 2$. Find the minimum size of a $1$-perfect bitrade in $H(n,q)$.
\end{problem}
For $q=2$ Problem \ref{ProblemPerfectBitrade} was essentially solved by Etzion and Vardy \cite{EtzionVardy} and Solov'eva \cite{Solov'eva}
(the results were formulated for more special cases of $1$-perfect bitrades embedded into perfect binary codes, but both proofs work in the general case).
In \cite{MogilnykhSolov'eva} Mogilnykh and Solov'eva for arbitrary $q\geq 2$ proved that the minimum size of a $1$-perfect bitrade in $H(q+1,q)$ is
$2\cdot q!$.

Now, using Theorem \ref{TheoremHamming}, we give a short solution of Problem \ref{ProblemPerfectBitrade} for $q=3$ and $q=4$.
Firstly, we need the following result.
\begin{lemma}[\cite{Valyuzhenich20}, Lemma 6]\label{f(T0,T1)}
Let $(T_0,T_1)$ be a $1$-perfect bitrade in a graph $G$. Then $f_{(T_0,T_1)}$ is a $(-1)$-eigenfunction of $G$.
\end{lemma}
Lemma \ref{f(T0,T1)} implies that we can consider Problem \ref{ProblemPerfectBitrade} only for $n=qm+1$, where $m\geq 1$ (because $-1$ is an eigenvalue of $H(n,q)$).

Suppose that $(T_0,T_1)$ is a $1$-perfect bitrade in $H(qm+1,q)$. By Lemma \ref{f(T0,T1)} we have that $f_{(T_0,T_1)}$ is a $(-1)$-eigenfunction of $H(qm+1,q)$. We note that $-1=\lambda_{(q-1)m+1}(qm+1,q)$.
Applying Theorem \ref{TheoremHamming} for $n=qm+1$ and $i=j=(q-1)m+1$, we obtain that
$$|S(f_{(T_0,T_1)})|\geq 2^{(q-1)m+1}\cdot (q-1)^{m}$$ for $q\geq 4$ and
$$|S(f_{(T_0,T_1)})|\geq 2^{m+1}\cdot 3^{m}$$ for $q=3$.
Consequently, we have
\begin{equation}\label{BitradeBoundOne}
|T_0|+|T_1|=|S(f_{(T_0,T_1)})|\geq 2^{(q-1)m+1}\cdot (q-1)^{m}
\end{equation}
for $q\geq 4$ and
\begin{equation}\label{BitradeBoundTwo}
|T_0|+|T_1|=|S(f_{(T_0,T_1)})|\geq 2^{m+1}\cdot 3^{m}
\end{equation}
for $q=3$.
On the other hand, in \cite{MogilnykhSolov'eva} Mogilnykh and Solov'eva for arbitrary $q\geq 2$ showed the existence of $1$-perfect bitrades in $H(qm+1,q)$ of size
$2\cdot (q!)^m$. Thus, the bounds (\ref{BitradeBoundOne}) and (\ref{BitradeBoundTwo}) are sharp for $q=4$ and $q=3$ respectively, and we obtain a solution of Problem \ref{ProblemPerfectBitrade} for $q\in \{3,4\}$. 
Finally, we note that Theorem \ref{TheoremHamming} implies that for $q\geq 5$ optimal $(-1)$-eigenfunctions of the Hamming graph $H(qm+1,q)$
do not correspond to its $1$-perfect bitrades (in this case we have a characterization of all optimal $(-1)$-eigenfunctions).
So, this approach does not work for $q\geq 5$.
\section{Doob graph}\label{SectionDoob}
In this section, we give a survey of results on MS-problem for the Doob graph.
The {\em Shrikhande graph} $\mathrm{Sh}$ is the Cayley graph on the group $\mathbb{Z}_{4}^{2}$  with the generating set
$\{\pm(0,1),\pm(1,0),\pm(1,1)\}$.
\begin{figure}[H]
\begin{center}
\includegraphics[scale=1.03]{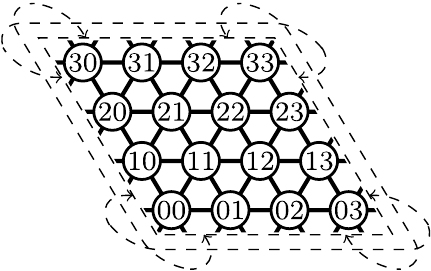}
\end{center}
\caption{The Shrikhande graph.}\label{Figure5}
\end{figure}
The {\em Doob graph} $D(m,n)$, where $m>0$, is the Cartesian product of $m$ copies of the Shrikhande graph and $n$ copies of the complete graph $K_4$.
In other words, we have $D(m,n)=\mathrm{Sh}^m\square K_4^n$.
This graph is a distance-regular graph with the same parameters as the Hamming graph $H(2m+n,4)$.
The Doob graph $D(m,n)$ has $2m+n+1$ distinct eigenvalues $\lambda_{i}(m,n)=6m+3n-4i$, where $0\leq i\leq 2m+n$.
In \cite{Bespalov} Bespalov proved that the minimum cardinality of the support of a $\lambda_{1}(m,n)$-eigenfunction of $D(m,n)$ is $6 \cdot 4^{2m+n-2}$ and obtained a characterization of optimal $\lambda_{1}(m,n)$-eigenfunctions. He also showed that the minimum cardinality of the support of a $\lambda_{2m+n}(m,n)$-eigenfunction of $D(m,n)$ is $2^{2m+n}$ and obtained a characterization of optimal $\lambda_{2m+n}(m,n)$-eigenfunctions.
In what follows, in this section we will consider the results obtained in \cite{Bespalov}.

Now, we discuss constructions of optimal $\lambda_{1}(m,n)$-eigenfunctions and $\lambda_{2m+n}(m,n)$-eigenfunctions.
It is interesting that as in the case of the Hamming graph such functions are constructed as a tensor product of several elementary optimal eigenfunctions.
Firstly, we define two sets of elementary optimal eigenfunctions.

\begin{figure}[H]
\begin{center}
\includegraphics[scale=1.03]{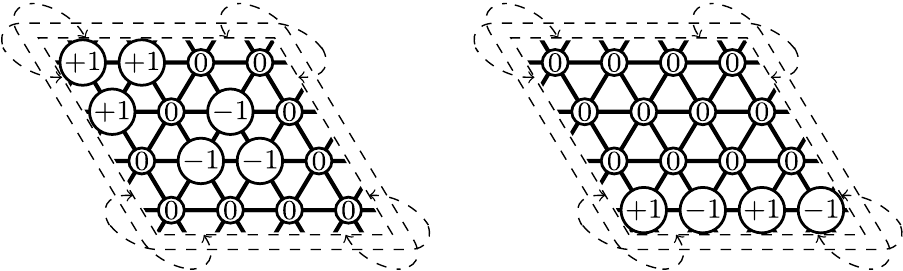}
\end{center}
\caption{Functions $p_{(0,3)}$ and $r_{(0,0),(0,1)}$.}\label{Figure6}
\end{figure}

For $a\in \mathbb{Z}_{4}^{2}$ we define the function $p_a$ on the vertices of the Shrikhande graph by the following rule:
$$
p_{a}(x)=\begin{cases}
1,&\text{if $x\in\{a+(3,1), a+(3,2), a+(2,1)\}$;}\\
-1,&\text{if $x\in \{a+(2,3),a+(1,2), a+(1,3)\}$;}\\
0,&\text{otherwise.}
\end{cases}
$$
We note that the support of $p_a$ consists of two disjoint copies of the complete graph $K_3$. The function $p_{(0,3)}$ is shown in Figure \ref{Figure6}. Denote $P=\{p_{a}~|~a\in{\mathbb{Z}_{4}^{2}}\}$.

For $a\in \mathbb{Z}_{4}^{2}$ and $b\in \{(0,1),(1,0),(1,1)\}$ we define the function $r_{a,b}$ on the vertices of the Shrikhande graph by the following rule:
$$
r_{a,b}(x)=\begin{cases}
1,&\text{if $x\in\{a, a+2b\}$;}\\
-1,&\text{if $x\in \{a+b, a+3b\}$;}\\
0,&\text{otherwise.}
\end{cases}
$$
We note that the vertices from the support of $r_{a,b}$ form a cycle of length $4$.  The function $r_{(0,0),(0,1)}$ is shown in Figure \ref{Figure6}. Denote $R=\{r_{a,b}~|~a\in{\mathbb{Z}_{4}^{2},b\in \{(0,1),(1,0),(1,1)\}}\}$.
We will also use the sets of functions $A_4$ and $C_4$ defined in Section \ref{SectionHamming}.

For $m>0$ denote by $I^{m,n}$ the function that is defined on the vertices of $D(m,n)$ and is identically equal to $1$.
For $n\geq 1$ denote by $I^n$ the function that is defined on the vertices of $H(n,4)$ and is identically equal to $1$.

Now, we define two classes of optimal $\lambda_{1}(m,n)$-eigenfunctions and one class of optimal $\lambda_{2m+n}(m,n)$-eigenfunctions.

We say that a function $f$ defined on the vertices of $D(m,n)$ belongs to the class $\mathcal{G}_{1}(m,n)$ if
$f=c\cdot g_1 \ldots g_{m}\cdot I^{n}$, where $c$ is a real non-zero constant, $g_i\in P$ for some $i \in \{1,\ldots,m\}$ and
$g_j=I^{1,0}$ for any $j\in \{1,\ldots,m\}\setminus i$.

Let $n\geq 2$. We say that a function $f$ defined on the vertices of $D(m,n)$ belongs to the class $\mathcal{G}_{2}(m,n)$ if
$f=c\cdot I^{m,0}\cdot h_1 \ldots h_{n-1}$, where $c$ is a real non-zero constant, $h_i\in A_4$ for some $i \in \{1,\ldots,n-1\}$ and
$h_j=I^{1}$ for any $j\in \{1,\ldots,n-1\}\setminus i$.

We say that a function $f$ defined on the vertices of $D(m,n)$ belongs to the class $\mathcal{G}_{3}(m,n)$ if
$f=c\cdot g_1\ldots g_m\cdot h_1\ldots h_n$, where $c$ is a real non-zero constant, $g_i\in R$ for any $i \in \{1,\ldots,m\}$ and $h_j\in C_4$ for any $j \in \{1,\ldots,n\}$.

The main results proved in \cite{Bespalov} are the following.
\begin{theorem}[\cite{Bespalov}, Theorem 1]\label{TheoremDoob1}
Let $f$ be a $\lambda_{1}(m,n)$-eigenfunction of $D(m,n)$, where $m>0$.
Then $|S(f)| \ge 6 \cdot 4^{2m+n-2}$. Moreover, if $|S(f)|=6 \cdot 4^{2m+n-2}$, then the following statements hold:
\begin{enumerate}
  \item If $n\geq 2$, then $f\in \mathcal{G}_1(m,n)$ or $f\in \mathcal{G}_2(m,n)$.
  \item If $n\in \{0,1\}$, then $f\in \mathcal{G}_1(m,n)$.
\end{enumerate}
\end{theorem}

\begin{theorem}[\cite{Bespalov}, Theorem 2]\label{TheoremDoob2}
Let $f$ be a $\lambda_{2m+n}(m,n)$-eigenfunction of $D(m,n)$, where $m>0$.
Then $|S(f)| \ge 2^{2m+n}$. Moreover, if $|S(f)|=2^{2m+n}$, then $f\in \mathcal{G}_3(m,n)$.
\end{theorem}
\begin{remark}
We note that the bound proved in Theorem \ref{TheoremDoob2} can also be obtained by applying the weight distribution bound for the smallest eigenvalue of the Doob graph.
\end{remark}
\section{Johnson graph}\label{SectionJohnson}
In this section, we give a survey of results on MS-problem for the Johnson graph.
The {\em Johnson graph} $J(n,\omega)$ is defined as follows. The vertices of $J(n,\omega)$ are the binary vectors of length $n$ with $\omega$ ones; and two vertices are adjacent if they have exactly $\omega-1$ common ones.
The Johnson graph $J(n,\omega)$ has $\omega+1$ distinct eigenvalues $\lambda_{i}(n,\omega)=(\omega-i)(n-\omega-i)-i$, where $0\leq i\leq \omega$.
In \cite{VMV} Vorob'ev et al. showed that for a fixed $\omega$ and sufficiently large $n$ the minimum cardinality of the support of a $\lambda_{i}(n,\omega)$-eigenfunction of $J(n,\omega)$ is $2^{i}\cdot \binom{n-2i}{\omega-i}$ and obtained a characterization of optimal $\lambda_{i}(n,\omega)$-eigenfunctions. Thus, MS-problem for the Johnson graph is asymptotically solved for all eigenvalues.

Now we discuss the main results obtained in \cite{VMV}. Firstly, we define the function $f^{i,\omega,n}$ on the vertices of the Johnson graph $J(n,\omega)$ by the following rule:
$$f^{i,\omega,n}(x_1,\ldots,x_n)=\begin{cases}
1,&\text{if $x_{2k-1}+x_{2k}=1$ for any $1\leq k\leq i$ and $x_1+x_3+\ldots+x_{2i-1}$ is even;}\\
-1,&\text{if $x_{2k-1}+x_{2k}=1$ for any $1\leq k\leq i$ and $x_1+x_3+\ldots+x_{2i-1}$ is odd;}\\
0,&\text{otherwise.}
\end{cases} $$
So, the support of $f^{i,\omega,n}$ consists of binary vectors $(x_1,\ldots,x_n)$ of weight $\omega$ such that the product $(x_1-x_2)\cdot \ldots \cdot (x_{2i-1}-x_{2i})$ is not equal to zero.
In \cite[Proposition 1]{VMV} it was shown that $f^{i,\omega,n}$ is a $\lambda_{i}(n,\omega)$-eigenfunction of $J(n,\omega)$ and
$|S(f^{i,\omega,n})|=2^{i}\cdot \binom{n-2i}{\omega-i}$. The main result proved in \cite{VMV} is the following.

\begin{theorem}[\cite{VMV}, Theorem 4]\label{TheoremJohnson}
Let $i$ and $\omega$ be positive integers, $\omega\ge i$. There is $n_0(i,\omega)$ such
that for all $n\geq n_0(i,\omega)$
 and any $\lambda_i(n,\omega)$-eigenfunction $f$ of
$J(n,\omega)$ the following holds:
\begin{equation}\label{BoundJohnson}
|S(f)|\geq 2^{i}\cdot \binom{n-2i}{\omega-i},
\end{equation}
and any function that attains the bound (\ref{BoundJohnson}) is equivalent to $f^{i,\omega,n}$ up to a permutation of coordinate positions and the multiplication by a scalar.
\end{theorem}
\begin{remark}
We note that the bound (\ref{BoundJohnson}) for $i=\omega$ and arbitrary $n$ can also be obtained by applying the weight distribution bound for the smallest eigenvalue of the Johnson graph.
\end{remark}
Now, we discuss the main ideas of the proof of Theorem \ref{TheoremJohnson}.

Let $f$ be a real-valued function defined on the vertices of the
Johnson graph $J(n,\omega)$ and let $j_1,j_2\in
\{1,2,\dots,n\}$, $j_1<j_2$. We define a function
$f_{j_1,j_2}$ on the vertices of $J(n-2,\omega-1)$ as follows: for any vertex $y=(y_1,y_2, \dots
,y_{j_1-1},y_{j_1+1}, \dots ,y_{j_2-1},y_{j_2+1},\dots,y_n)$ of
$J(n-2,\omega-1)$
$$
\begin{gathered}
f_{j_1,j_2}(y)=f(y_1,y_2, \dots,y_{j_1-1},1,y_{j_1+1}, \dots ,y_{j_2-1},0,y_{j_2+1},\dots,y_n)-\\
-f(y_1,y_2, \dots ,y_{j_1-1},0,y_{j_1+1}, \dots ,y_{j_2-1},1,y_{j_2+1}, \dots ,y_n).
\end{gathered}
$$
One of the important ingredients in the proof of Theorem \ref{TheoremJohnson} is the following.
\begin{lemma}[\cite{VMV}, Lemma 1]\label{ReductionJohnson}
Let f be a $\lambda_i(n,\omega)$-eigenfunction of $J(n,\omega)$, where  $j_1,j_2\in\{1,2,\dots,n\}$ and $j_1<j_2$.
Then $f_{j_1,j_2}$ is a
$\lambda_{i-1}(n-2,\omega-1)$-eigenfunction of $J(n-2,\omega-1)$ or the all-zero function.
\end{lemma}
Lemma \ref{ReductionJohnson} is a very useful tool for studying of eigenfunctions of the Johnson graph.
In particular, this lemma allows to apply induction on $n$, $\omega$ and $i$ (we can use the induction assumption for the function $f_{j_1,j_2}$).
Moreover, we suppose that Lemma \ref{ReductionJohnson} can be useful not only for the MS-problem but also for other problems.
For example, recently in \cite{Vorob'ev} Vorob'ev applied Lemma \ref{ReductionJohnson} for characterization of equitable $2$-partitions of the Johnson graph with the eigenvalue $\lambda_2(n,\omega)$.
Finally, we note that Lemma \ref{ReductionJohnson} is an analogue of Lemma \ref{ReductionHamming} (see Subsection \ref{SubSectionResults}).

Let $v=(v_1,\ldots,v_n)$ be a real non-zero vector such that $v_1+\ldots+v_n=0$. We define the function $f^{v}$ on the vertices of the Johnson graph $J(n,\omega)$ by the following rule:
$$f^{v}(x_1,\ldots,x_n)=\sum_{1\leq i\leq n~:~x_i=1}v_i.$$
For $k\in\{1,\dots,n-1\}$ denote
$$v^k=(\underbrace{1,\ldots,1}_{k},\underbrace{\frac{-k}{n-k},\ldots,\frac{-k}{n-k}}_{n-k}).$$
In \cite{MVV} Mogilnykh  et al. proved the following improvement of Theorem \ref{TheoremJohnson} for $\lambda_{1}(n,\omega)$-eigenfunctions of $J(n,\omega)$.
\begin{theorem}[\cite{MVV}, Theorem 1]
Let $f$ be an optimal $\lambda_{1}(n,\omega)$-eigenfunction of $J(n,\omega)$, where $n\geq 2\omega$ and $\omega\geq 2$.
Then $f$ is $f^{1,\omega,n}$ or $f^{v^k}$ for some $k\in\{2,\dots,n-2\}$ such that $\frac{k\omega}{n}\in \mathbb{N}$ up to a permutation of coordinate positions and the multiplication
by a scalar.
\end{theorem}

\section{Grassmann Graph}\label{sec:grassmann}

In this section, we give a survey of results on MS-problem for the Grassmann graph. The Grassmann graph $J_q(N,m)$ is a distance-regular graph with the vertex set consisting of all $m$-dimensional subspaces of a vector space of dimension $N$ over a finite field $\mathbb{F}_q$. Two vertices are adjacent whenever the corresponding subspaces intersect in a $(m-1)$-dimensional subspace. 

MS-problem for the minimum eigenvalue of the Grassmann graph was studied in \cite{KrotovMogilnykhPotapov}. But it is interesting that this problem can be tracked earlier to the works \cite{J84, Cho99, Cho98} where it was considered in terms of finding the minimum null $t$-designs of the lattices of subspaces over a finite field. 
In \cite{J84} G. D. James made a conjecture about the minimum support size of non-zero null $t$-designs of the lattices of subspaces over a finite field. S. Cho confirms the conjecture in \cite{Cho99} and in \cite{Cho98} characterizes all the null $t$-designs with minimum supports in terms of maximal isotropic spaces of some bilinear form. 


Coming back to the Grassmann graph, we obtain the following theorem that gives us the characterization of optimal $\lambda_D$-eigenfunctions for the Grassmann graph (compare with Theorems~1,2 from \cite{Cho98} and Theorem~5 from \cite{KrotovMogilnykhPotapov}). For more details about null $t$-designs and totally isotropic spaces the reader is referred to \cite{Cho98} and Chapter 18 of \cite{J84}.

\begin{theorem}\label{th:grassmann_char}
	Suppose $f$ is an optimal $\lambda_D$-eigenfunction of the Grassmann graph $J_q(N,m)$, where $N\ge 2m$ and $\lambda_D$ is its minimum eigenvalue. Then the cardinality of its support is $\sum\limits_{i=0}^D\begin{bmatrix}m\\i\end{bmatrix}_{q}\cdot q^{i(i-1)/2}$ which is also equal to the value of the weight distribution bound and the non-zeros of the function $f$ correspond to the maximal totally isotropic subspaces of a $2m$-dimensional space, equipped with a bilinear form $B$ with a Gram matrix $\begin{pmatrix}
	\mathbf{0} & E_m\\
	E_m & \mathbf{0}
	\end{pmatrix}$ up to the equivalence (or, equivalently, with respect to a non-degenerate quadratic form $Q$).
\end{theorem}

Thus for the minimum eigenvalue of the Grassmann graph $J_q(N,m)$ MS-problem is solved and the weight distribution bound is achieved.

\section{Bilinear Forms Graph}\label{sec:bfg}

In this section, we give a survey of results on MS-problem for bilinear forms graph. More details can be found \cite{SotnikovaBilinear}. The bilinear forms graph $\mathrm{Bil}_q(n,m)$ is a distance-regular graph with the vertex set $V$ consisting of all $n\times m$ matrices over a finite field $\mathbb{F}_q$ and two vertices being adjacent when their matrix difference has a rank $1$. For the sake of convenience, we will further suppose that $m\le n$. Thus the diameter $D$ of the bilinear forms graph $\mathrm{Bil}_q(n,m)$ is equal to $m$.

Here as well as in the previous section we consider MS-problem only for the case of minimum eigenvalue $\lambda_D$. In this case we have the following lower bound for the minimum support cardinality: \[\sum\limits_{i=0}^{m} \begin{bmatrix}m\\i\end{bmatrix}_{q}\cdot q^{i(i-1)/2}\]
It is interesting that the weight distribution for bilinear forms graph coincides with that of the Grassmann graph. Later we will see the importance of this connection.

The key idea here is that bilinear forms graph belongs to a family of so-called Delsarte cliques graphs (each edge lies in a constant number of Delsarte cliques). Recall that a clique in a distance-regular graph of degree $k$ is called Delsarte clique if it consists of exactly $1-k/\lambda_D$ vertices. For more details about Delsarte cliques graphs, the reader is referred to \cite{BHK07}. 

This property leads to the following observations:
\begin{itemize}
    \item Theorem~2 from \cite{KrotovMogilnykhPotapov} implies that for a Delsarte cliques graph $G$ a function $f$ is a $\lambda_D$-eigenfunction of $G$ if and only if for every Delsarte clique $C$ it holds $\sum\limits_{v\in C}f(v)=0$.
    
    \item Theorem~3 from \cite{KrotovMogilnykhPotapov} tells us that for a Delsarte clique graph $G$ in case of $D=2$ if a weight distribution bound is achieved then non-zeros of optimal $\lambda_D$-eigenfunction induce a complete bipartite graph. Note that for bilinear forms graph $\mathrm{Bil}_q(2,2)$ we have $\lambda_D=-q-1$, thus non-zeros of optimal $\lambda_D$-eigenfunction achieving weight distribution bound induce a complete bipartite graph $K_{q+1, q+1}$ if such a function exists.
\end{itemize}

It appears that in case of strongly regular bilinear forms graphs $\mathrm{Bil}_p(2,\,2)$ (those with $D=2$) the weight distribution bound can be achieved. An explicit construction of an optimal $\lambda_D$-eigenfunction can be found in \cite{SotnikovaBilinear}. Below are the statements that summarize this construction, but first let us introduce additional notation. Suppose $a_1$ is a generating element of the multiplicative group $\mathbb{F}_p^*$. Denote 

\[
a_0=0,\quad a_2=a_1^2,\quad \ldots,\quad a_{p-2}=a_1^{p-2}, \quad a_{p-1}=a_1^{p-1}=1
\]
\[
e_*=[0,1],\quad e_0=[1,0], \quad e_1=[1,a_1],\quad \ldots, \quad e_{p-1}=[1,a_{p-1}]
\]

\begin{theorem}[\cite{SotnikovaBilinear}, Theorem~3]\label{bfg:minsup_exists}
	Let $\mathrm{Bil}_p(2,2)$ be a bilinear forms graph over a prime field $\mathbb{F}_p$. For any $\nu\in \mathbb{F}_p$, such that $\nu\ne -\xi^2$ for all $\xi \in \mathbb{F}_p$, and $b_i=\frac{1}{a_i^2\nu+1}$ the independent set
	\[\mathcal{N}=\Bigl\{\begin{bmatrix}0\\1\end{bmatrix}e_*,\; b_0\begin{bmatrix}1\\a_0\nu\end{bmatrix}e_0,\; \ldots,\; b_{p-1}\begin{bmatrix}1\\a_{p-1}\nu\end{bmatrix}e_{p-1} \Bigr\}
	\]
	together with
	\[\mathcal{P}=\Bigl\{\begin{bmatrix}0\\0\end{bmatrix}e_*;\; b_0\begin{bmatrix}1\\a_0\nu\end{bmatrix}e_0+b_0\begin{bmatrix}-a_0\\1\end{bmatrix}e_*;\; \ldots;\; b_{p-1}\begin{bmatrix}1\\a_{p-1}\nu\end{bmatrix}e_{0}+b_{p-1}\begin{bmatrix}-a_{p-1}\\1\end{bmatrix}e_* \Bigr\}
	\]
	form non-zeros of $\lambda_D$-eigenfunction $f$ as two parts of a complete bipartite graph $K_{p+1,\,p+1}$ and 
	\[
	f(v)=\begin{cases}
	c, & \text{for }v\in \mathcal{P}, \\
	-c, & \text{for }v\in \mathcal{N},\\
	0, & \textit{else}
	\end{cases}
	\]
	for some constant $c\ne 0$.
\end{theorem}

Let us illustrate this theorem with some small example. Consider a bilinear forms graph $\mathrm{Bil}_3(2,2)$. Using the construction above we obtain the following sets:

\[
\begin{gathered}
     \mathcal{N}=\Bigl\{
    \begin{bmatrix}0 & 0\\0 & 1\end{bmatrix},
    \begin{bmatrix}1 & 0\\0 & 0\end{bmatrix},
    \begin{bmatrix}2 & 1\\1 & 2\end{bmatrix},
    \begin{bmatrix}2 & 2\\2 & 2\end{bmatrix}
    \Bigr\},\\
     \mathcal{P}=\Bigl\{
    \begin{bmatrix}0 & 0\\0 & 0\end{bmatrix},
    \begin{bmatrix}1 & 0\\0 & 1\end{bmatrix},
    \begin{bmatrix}2 & 2\\1 & 2\end{bmatrix},
    \begin{bmatrix}2 & 1\\2 & 2\end{bmatrix}
    \Bigr\}.
\end{gathered}
\]

Here under the notation above $a_0=0$, $a_1=2$, $a_2=1$, $e_*=[0,1]$, $e_0=[1,0]$, $e_1=[1,2]$, $e_2=[1,1]$, $\nu=1$, $b_0=1$, $b_1=2$, $b_2=2$.

Thus we proved that there exists a family of optimal $\lambda_D$-eigenfunctions of the bilinear forms graph $\mathrm{Bil}_p(2,\,2)$ over a prime field $\mathbb{F}_p$ that achieve the lower bound. However the construction described above does not provide the full characterization of all optimal $\lambda_D$-eigenfunctions.

What happens if we look at bilinear forms graphs of larger diameter? It appears that the weight distribution bound cannot be achieved. And for proving this the connection between bilinear graphs and the Grasssmann graphs comes in handy. The bilinear forms graph $\mathrm{Bil}_q(n,m)$ with $m\le n$ can be considered as a subgraph of the Grassman graph $J_q(n+m, m)$ as follows: given a fixed subspace $W$ of dimension $n$, all $m$-spaces $U$ such that $U\cap W=0$ are the vertices of $\mathrm{Bil}_q(n,m)$. This embedding leads to the following result about the Delsarte cliques of these graphs (see Lemma~8 from \cite{SotnikovaBilinear}):

\begin{lemma}\label{l:bfg_grassmann_delsarte_cliques}
	Delsarte cliques of bilinear forms graph $\mathrm{Bil}_q(n,m)$ are embedded in Delsarte cliques of a Grassmann graph $J_q(n+m,m)$ in the sense that for any Delsarte cliques $C$ and $\widehat{C}$ of a bilinear forms graph and the Grassmann graph correspondingly, either $C\subset \widehat{C}$ or $C\cap \widehat{C}=\emptyset$. 
\end{lemma}

Since for any $\lambda_D$-eigenfunction the sum of its values over a Delsarte clique is zero, from the previous Lemma we immediately obtain the following Corollary which simply tells us that we can extend eigenfunctions of bilinear forms graph to those of the Grassmann graph:

\begin{corollary}\label{c:bil_grassmann_eigens}
	Suppose $f$ is a $\lambda_D$-eigenfunction of a bilinear forms graph $\mathrm{Bil}_q(n,m)$. Then $\widehat{f}$ is an eigenfunction of the Grassmann graph $J_q(n+m,m)$, where 
	\[
	\widehat{f}(\overline{M})  = 
	\begin{cases}
	f(M), & \text{if }M\in V(\mathrm{Bil}_q(n,m)) \\
	0, & \text{else}
	\end{cases}
	\]
\end{corollary}

This corollary is crucial for the final result:

\begin{theorem}[\cite{SotnikovaBilinear}, Theorem~7]\label{th:bfg_not_exist}
	Let $\mathrm{Bil}_q(n,m)$ be a bilinear forms graph of diameter $D\ge 3$. Then the minimum support of an eigenfunction corresponding to the minimum eigenvalue does not achieve the weight distribution bound.
\end{theorem}

The main idea behind the proof of this theorem can be described as follows. Suppose the opposite holds and $f$ is an optimal $\lambda_D$-eigenfunction that achieves the weight distribution bound. Under the notation of Corollary~\ref{c:bil_grassmann_eigens}, $\widehat{f}$ is an optimal $\lambda_D$-eigenfunction of the Grassmann graph $J_q(n+m,\,m)$. According to the Theorem~\ref{th:grassmann_char} characterizing optimal eigenfunctions of the Grassmann graphs, the non-zeros of $\widehat{f}$ correspond to the maximal totally isotropic spaces of a non-degenerate quadratic form $Q$. Now we recall the graphs embedding construction: there exists a subspace $W$ of dimension $n$ that trivially intersects with all the maximal totally isotropic subspaces. A well-known corollary from the Chevalley theorem states that  any non-degenerate quadratic form is isotropic on a vector space of dimension not less that $3$ over the finite field $\mathbb{F}_q$ (here the diameter of a graph plays its role). Thus there exists a non-zero vector $w\in W$ such that $Q(w)=0$, therefore ${<}w{>}$ is a $1$-dimensional totally isotropic space and, hence, is contained in a maximal totally isotropic subspace. This contradicts the trivial intersection of $W$ with all the maximal totally isotropic subspaces.

According to this theorem optimal $\lambda_D$-eigenfunctions of $\mathrm{Bil}_q(n,\,m)$ do not satisfy the weight distribution bound. This lead to an open MS-problem for bilinear forms graphs of diameter $D\ge 3$. 

\section{Paley graph}\label{SectionPaley}
In this section, we give a survey of results on MS-problem for the Paley graph.
Let $q$ be an odd prime power, where $q\equiv 1(4)$.
The Paley graph $P(q)$ is the Cayley graph on the additive group $\mathbb{F}_{q}^+$ of the finite field $\mathbb{F}_q$ with the generating set of all squares in the multiplicative group $\mathbb{F}_{q}^*$.
This graph is a strongly regular with parameters $(q,\frac{q-1}{2},\frac{q-5}{4},\frac{q-1}{4})$.
The eigenvalues of $P(q)$ are $\lambda_0=\frac{q-1}{2}$, $\lambda_1=\frac{-1+\sqrt{q}}{2}$ and $\lambda_2=\frac{-1-\sqrt{q}}{2}$.
In \cite{GoryainovKabanovShalaginovValyuzhenich} Goryainov et al. for $i\in \{1,2\}$ proved that the minimum cardinality of the support of a $\lambda_i$-eigenfunction of $P(q^2)$, where $q$ is an odd prime power, is $q+1$.
In what follows, in this section we will discuss the results obtained in \cite{GoryainovKabanovShalaginovValyuzhenich}.

Let $q$ be an odd prime power and let $\beta$ be a primitive element of the finite field  $\mathbb{F}_{q^2}$.
Denote $\omega=\beta^{q-1}$, $Q_0=\langle \omega^2 \rangle$ and  $Q_1=\omega \langle \omega^2 \rangle$.
We define the function $f_{\beta}$ on the vertices of the Paley graph $P(q^2)$ by the following rule:
$$
f_{\beta}(x)=\begin{cases}
1,&\text{if $x\in Q_0$;}\\
-1,&\text{if $x\in Q_1$;}\\
0,&\text{otherwise.}
\end{cases}
$$
One of the main results proved in \cite{GoryainovKabanovShalaginovValyuzhenich} is the following.
\begin{theorem}[\cite{GoryainovKabanovShalaginovValyuzhenich}, Theorem 2]\label{TheoremPaley}
Let $q$ be an odd prime power and let $\beta$ be a primitive element of the finite field  $\mathbb{F}_{q^2}$.
Then the following statements hold:
\begin{enumerate}
  \item If $q\equiv 1(4)$, then $f_\beta$ is a $\lambda_2$-eigenfunction of $P(q^2)$ and $|S(f_\beta)|=q+1$.

  \item If $q\equiv 3(4)$, then $f_\beta$ is a $\lambda_1$-eigenfunction of $P(q^2)$ and $|S(f_\beta)|=q+1$.
\end{enumerate}
\end{theorem}
Since the Paley graph $P(q^2)$ is self-complementary, Theorem \ref{TheoremPaley} implies that for any  $i\in \{1,2\}$  $P(q^2)$
has $\lambda_i$-eigenfunction $f$ such that $|S(f)|=q+1$.
On the other hand, by the weight distribution bound we obtain that a $\lambda_2$-eigenfunction of $P(q^2)$ has at least $q+1$ non-zero values.
Since $P(q^2)$ is self-complementary, the same bound holds for a $\lambda_1$-eigenfunction of $P(q^2)$.
Thus, the minimum cardinality of the support of a $\lambda_i$-eigenfunction of $P(q^2)$, where $i\in \{1,2\}$, is $q+1$.

Now we discuss one interesting connection between the sets $Q_0$ and $Q_1$ and maximal cliques of the Paley graph $P(q^2)$.
The maximum possible size of a clique of $P(q^2)$ is $q$ (all cliques of such size are Delsarte cliques).
Blokhuis \cite{Blokhius84} determined all cliques and all cocliques of size $q$ in $P(q^2)$ and showed that
they are affine images of the subfield $\mathbb{F}_q$.
Baker et al. \cite{BakerEHW96} found maximal cliques of order $\frac{q+1}{2}$ and $\frac{q+3}{2}$ for $q \equiv 1(4)$ and $q \equiv 3(4)$ respectively,
but these cliques are not the only cliques of such size. Moreover, there are no known maximal cliques whose size
belongs to the gap from $\frac{q+1}{2}$ (from $\frac{q+3}{2}$, respectively) to $q$.
Kiermaier and Kurz \cite{KiermaierK} studied maximal integral point sets in affine planes over finite fields and found maximal cliques of size 
$\frac{q+3}{2}$ in $P(q^2)$ for $q \equiv 3(4)$.
Using the sets $Q_0$ and $Q_1$ defined above, Goryainov et al. \cite{GoryainovKabanovShalaginovValyuzhenich} constructed new maximal cliques of size $\frac{q+1}{2}$ and $\frac{q+3}{2}$ for $q \equiv 1(4)$ and $q \equiv 3(4)$ respectively in $P(q^2)$.
\begin{theorem}[\cite{GoryainovKabanovShalaginovValyuzhenich}, Theorem 1]\label{TheoremPaleyCliques}
Let $q$ be an odd prime power and let $\beta$ be a primitive element of the finite field  $\mathbb{F}_{q^2}$.
Then the following statements hold:
\begin{enumerate}
  \item If $q\equiv 1(4)$, then $Q_0$ and $Q_1$ are maximal cocliques of size $\frac{q+1}{2}$ in the graph $P(q^2)$.

  \item If $q\equiv 3(4)$, then $Q_0\cup \{0\}$ and $Q_1\cup \{0\}$ are maximal cliques of size $\frac{q+3}{2}$ in the graph $P(q^2)$.
\end{enumerate}
\end{theorem}
\section{The Star graph}\label{SectionStar}
In this section, we give a survey of results on MS-problem for the Star graph.
The {\em Star graph} $S_n$, $n\ge 3$, is the Cayley graph on the symmetric group $\rm{Sym}_n$  with the generating set 
$\{(1~i) ~|~ i \in \{2,\ldots,n\}\}$.
This graph is not distance-regular. The spectrum of the Star graph is integral \cite{ChapuyFeray,KrakovskiMohar}.
For $n\geq 4$, the eigenvalues of $S_n$ are $\pm (n-k)$, where $1\le k \le n$;
and the eigenvalues of $S_3$ are $\{-2,-1,1,2\}$.
The multiplicities of eigenvalues of the Star graph were studied in \cite{AKK16,KK15,K18}. In particular, explicit formulas for calculating multiplicities of eigenvalues $\pm(n-k)$, where $2\leq k \leq 12$, were found.
In \cite{KabanovKonstantinovaShalaginovValyuzhenich} Kabanov et al. found the minimum cardinality of the support of an $(n-2)$-eigenfunction of $S_n$ and obtained a characterization of optimal $(n-2)$-eigenfunctions for $n\ge 8$ and $n=3$.
In what follows, in this section we will consider the results obtained in \cite{KabanovKonstantinovaShalaginovValyuzhenich}.

Now, we discuss one construction of optimal $(n-2)$-eigenfunctions of the Star graph.
Let $i\in\{1,\ldots,n\}$ and $j,k\in\{2,\ldots,n\}$, where $j\neq k$.
We define the function $f_{i}^{j,k}$ on the vertices of the Star graph $S_n$ by the following rule:
$$
f_{i}^{j,k}(\pi)=\begin{cases}
1,&\text{if $\pi(j)=i$;}\\
-1,&\text{if $\pi(k)=i$;}\\
0,&\text{otherwise.}
\end{cases}
$$
In \cite[Lemma 2]{KabanovKonstantinovaShalaginovValyuzhenich} it was shown that $f_{i}^{j,k}$ is an $(n-2)$-eigenfunction of $S_n$ and  $|S(f_{i}^{j,k})|=2(n-1)!$.
Denote $$\mathcal{F}=\{f_{i}^{j,k}~|~i\in\{1,\ldots,n\},j,k\in\{2,\ldots,n\},j\neq k\}.$$
The main result proved in \cite{KabanovKonstantinovaShalaginovValyuzhenich} is the following.
\begin{theorem}[\cite{KabanovKonstantinovaShalaginovValyuzhenich}, Theorem 20]\label{TheoremStar}
Let $f$ be an $(n-2)$-eigenfunction of $S_n$, where $n\ge 8$ or $n=3$. Then $|S(f)|\ge 2(n-1)!$. Moreover, $|S(f)|=2(n-1)!$ if and only if $f=c\cdot \tilde{f}$, where $c$ is a real non-zero constant and $\tilde{f}\in \mathcal{F}$.
\end{theorem}
Now, we discuss the main ideas of the proof of Theorem \ref{TheoremStar}. Firstly, we need some definitions.

Let  $M=(m_{i,j})$ be a real $n\times n$ matrix.
We say that $M$ is {\em special} if $M$ is non-zero and the following conditions hold:
\begin{enumerate}
  \item $m_{i,1}=0$ for any $i\in{\{1,\ldots,n\}}$.
  \item $m_{1,j}=0$ for any $j\in{\{1,\ldots,n\}}$.
  \item $\sum_{j=1}^{n}m_{i,j}=0$ for any $i\in{\{1,\ldots,n\}}$.
\end{enumerate}
For a real $n\times n$ matrix $M=(m_{i,j})$ denote
$$g_M(n)=|\{\pi\in \rm{Sym}_n~|~\sum_{i=1}^{n}m_{i,\pi(i)}\neq 0\}|.$$
The key point of the proof of Theorem \ref{TheoremStar} is the following.
For an arbitrary $(n-2)$-eigenfunction $f$ of $S_n$ we can construct some special $n\times n$ matrix $M(f)$ and match the permutations from $\rm{Sym}_n$ with diagonals of $M(f)$ in such a way that the value of $f$ on a permutation $\pi$ is the sum of elements of the corresponding diagonal of $M(f)$.
In other words, we have the equality
\begin{equation}\label{EqStar}
|S(f)|=g_{M(f)}(n)
\end{equation}
for any $(n-2)$-eigenfunction $f$ of $S_n$.
This observation allows us to reduce MS-problem for the Star graph $S_n$ and its eigenvalue $n-2$ to the following extremal
problem on the set of all special $n\times n$ matrices.
\begin{problem}\label{MinSpecMatr}
Given a positive integer $n$, to find the minimum value of $g_M(n)$ for the class of special $n \times n$ matrices $M$.
\end{problem}
In \cite[Theorem 19]{KabanovKonstantinovaShalaginovValyuzhenich} for $n\ge 8$ and $n=3$ it was proved that $g_M(n)\ge 2(n-1)!$ for any special $n \times n$ matrix $M$. Moreover, in \cite[Theorem 19]{KabanovKonstantinovaShalaginovValyuzhenich} a classification of special matrices in the equality case was obtained. Using these results and the equality (\ref{EqStar}), we can finish the proof of Theorem \ref{TheoremStar}.

\section{Some remarks on optimal eigenfunctions of graphs}\label{SectionRemarks}
In this section, we give some observations on optimal eigenfunctions of graphs.

\begin{figure}[H]
\begin{center}
\includegraphics[scale=0.4]{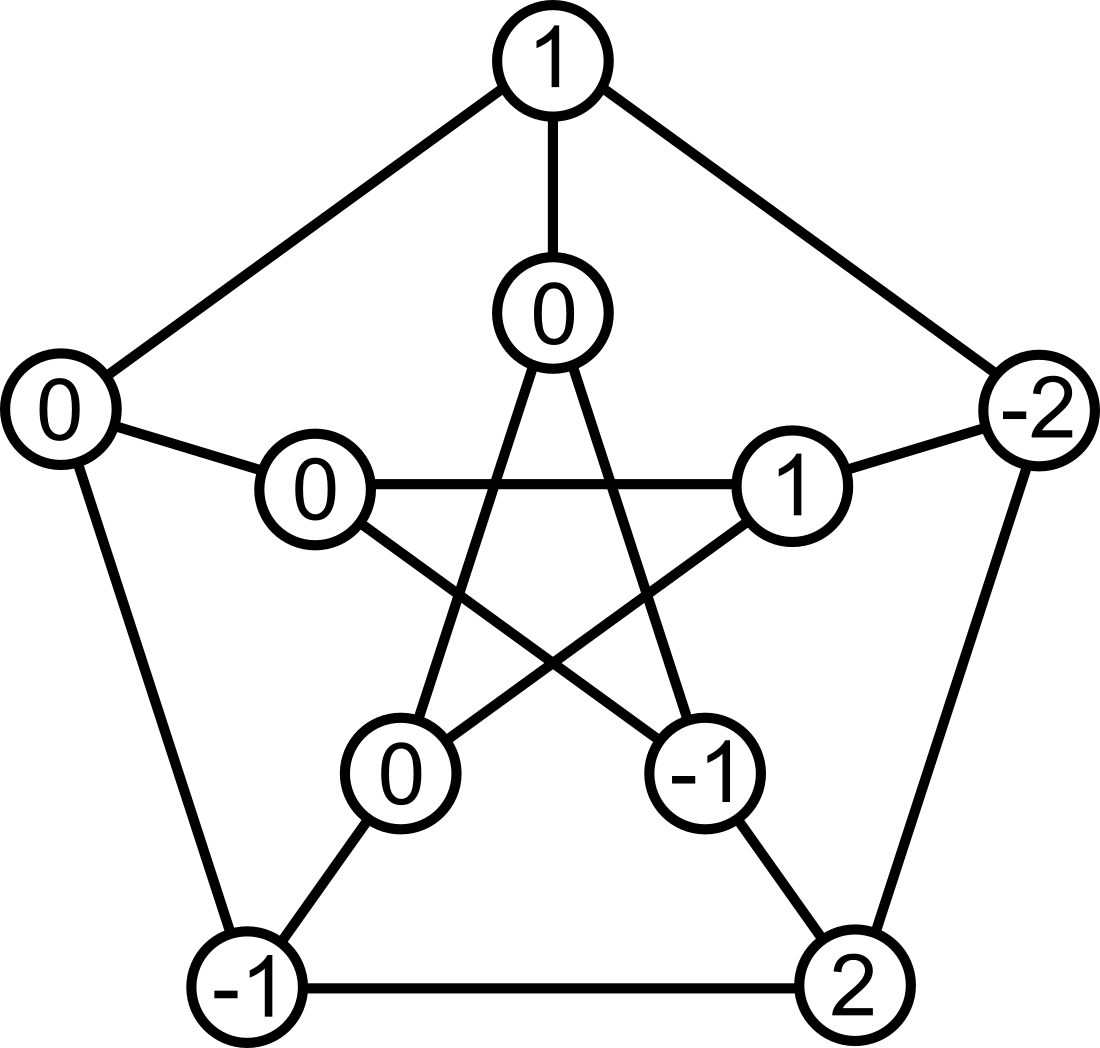}
\end{center}
\caption{Optimal $(-2)$-eigenfunction of the Petersen graph.}\label{Figure7}
\end{figure}

Recall that MS-problem is formulated for arbitrary real-valued functions from the corresponding eigenspace.
Surprisingly, in many cases optimal eigenfunctions take only three distinct values (for example, see Theorems
\ref{TheoremHamming}, \ref{TheoremDoob1}, \ref{TheoremDoob2}, \ref{TheoremJohnson}, \ref{TheoremStar}).
But, in general case it is not true.
For example, there are optimal $(-2)$-eigenfunctions of the Petersen graph that take five distinct values (see Figure \ref{Figure7}).

There is an interesting connection between optimal eigenfunctions corresponding to
the second largest eigenvalue of a given graph and completely regular codes in this graph.
In particular, an arbitrary optimal $\lambda_{1}(n,q)$-eigenfunction ($\lambda_{1}(n,\omega)$-eigenfunction) of the Hamming graph $H(n,q)$ (the Johnson graph $J(n,\omega)$) is the difference of the characteristic functions of two completely regular codes of covering radius $1$ (see \cite[Theorem 3]{Valyuzhenich} and \cite[Theorem 4]{VMV}).
The Star graph $S_n$ does not have completely regular codes of covering radius $1$ with the eigenvalue $n-2$. However, an arbitrary optimal $(n-2)$-eigenfunction of $S_n$ is the difference of the characteristic functions of two completely regular codes of covering radius $2$
(see \cite[Lemma 22]{KabanovKonstantinovaShalaginovValyuzhenich}).

\section{Open problems}\label{SectionOpenProblems}
In this section, we briefly recall the main results on MS-problem and formulate several open problems.

Recall that Problem \ref{GMSHamming} is completely solved for all $n\geq 1$ and $q\ge 2$. In particular, MS-problem for the Hamming graph $H(n,q)$ is solved for all eigenvalues. Moreover, a characterization of functions that are optimal in the space $U_{[i,j]}(n,q)$ was obtained for $q\ge 3$, $i+j\le n$ and $q\ge 5$, $i=j$, $i>\frac{n}{2}$. Taking into account these results, we formulate the following two problems for the Hamming graph.
\begin{problem}
Characterize functions that are optimal in the space $U_{[i,j]}(n,q)$ for the cases $q=2$ and $q\geq 3$, $i+j>n$ (in this problem we assume that $i<j$).
\end{problem}

\begin{problem}
Characterize optimal $\lambda_{i}(n,q)$-eigenfunctions of the Hamming graph $H(n,q)$ for $q\in \{3,4\}$ and $i>\frac{n}{2}$.
\end{problem}

MS-problem for the Doob graph $D(m,n)$ is solved for the second largest eigenvalue $\lambda_{1}(m,n)$ and the smallest eigenvalue $\lambda_{2m+n}(m,n)$.
So, it seems very interesting to consider the following question.
\begin{problem}
Solve MS-problem for the third largest eigenvalue $\lambda_{2}(m,n)$ of the Doob graph $D(m,n)$.
\end{problem}

MS-problem for the bilinear forms graph $\mathrm{Bil}_q(n,m)$  is solved for the smallest eigenvalue $\lambda_D$ in case $n=m=2$ and $q$ is prime. For bilinear forms graphs of larger diameters over the arbitrary field it is proved that the weight distribution bound cannot be attained. This leads to the following interesting questions:

\begin{problem}
For the bilinear forms graph $\mathrm{Bil}_q(n,m)$ of diameter $D$:
\begin{itemize}
    \item Characterize optimal $\lambda_D$-eigenfunctions in case of $D=2$ for a prime $q$ (including the case of $n\ne m$).
    \item Solve MS-problem for the smallest eigenvalue $\lambda_D$ in case of $D=2$ and arbitrary $q$.
    \item Solve MS-problem for the smallest eigenvalue $\lambda_D$ in case of $D\ge 3$ and arbitrary $q$.
\end{itemize}
\end{problem}

MS-problem for the Grassmann graph $J_q(N,m)$ is solved for the smallest eigenvalue $\lambda_D$. Since the Grassmann graph can be considered as a $q$-analogue of the Johnson graph it may be interesting to consider the following question:

\begin{problem}
Solve MS-problem for the second largest eigenvalue of the Grassmann graph $J_q(N,m)$.
\end{problem}

MS-problem for the Paley graph $P(q^2)$ is solved for both non-principal eigenvalues. We formulate the following problem for optimal eigenfunctions.
\begin{problem}
Characterize optimal $\lambda_1$-eigenfunctions and $\lambda_2$-eigenfunctions of the Paley graph $P(q^2)$.
\end{problem}

MS-problem for the Star graph $S_n$ is solved only for the second largest eigenvalue. So, the following question is very natural.
\begin{problem}
Solve MS-problem for the third largest eigenvalue of the Star graph $S_n$.
\end{problem}

At the end of this section we also would like to bring the attention of the reader to the following problems:
\begin{problem}
For distance-regular graphs find the conditions for the weight distribution bound to be achieved.
\end{problem}

\begin{problem}
For distance-regular graphs find a sharper lower bound on the cardinality of a graph eigenfunction support than the weight distribution bound.
\end{problem}

\begin{problem}
Find a lower bound on the cardinality of a graph eigenfunction support for the Cayley graphs.
\end{problem}

\section{Acknowledgements}
The authors are grateful  to Evgeny Bespalov, Denis Krotov, Vladimir Potapov and Konstantin Vorob'ev for helpful discussions.

\end{document}